\documentclass[11pt]{article} 

\usepackage{IS_via_SOC_manuscript}
\usepackage[section]{placeins}
\usepackage{tabularx,ragged2e,booktabs,caption}
\usepackage{hyperref}
\usepackage{algorithm}
\usepackage{footnote}
\numberwithin{equation}{section}
\numberwithin{figure}{section}
\numberwithin{table}{section}
\numberwithin{algorithm}{section}
\usepackage{authblk} 

\newcommand{\rset}{\mathbb{R}}
\newcommand{\nset}{\mathbb{N}}
\newcommand{\zset}{\mathbb{Z}}

\newcommand{\PERIOD}{.}
\newcommand{\COMMA}{,}

\newcommand{\Ordo}[1]{{\mathcal{O}}\left(#1\right)}
\newcommand{\ordo}[1]{{o}\left(#1\right)}

\makeatletter
\def\BState{\State\hskip-\ALG@thistlm}
\makeatother

\title{Learning-Based Importance Sampling via Stochastic Optimal Control for Stochastic Reaction Networks}

\author[1]{Chiheb Ben Hammouda}
\author[2]{Nadhir Ben Rached}
\author[3,4]{Ra\'ul Tempone}
\author[1]{Sophia Wiechert\thanks{wiechert@uq.rwth-aachen.de}}
\affil[1]{RWTH Aachen University, Chair of Mathematics for Uncertainty Quantification, Aachen, Germany.}
\affil[2]{University of Leeds, School of Mathematics, Leeds, UK.}
\affil[3]{King Abdullah University of Science and Technology (KAUST), Computer, Electrical and Mathematical Sciences \& Engineering Division (CEMSE), Thuwal, Saudi Arabia}
\affil[4]{RWTH Aachen University, Alexander von Humboldt Professor in Mathematics for Uncertainty Quantification, Aachen, Germany.}

\begin{document} 
    \date{}
\maketitle

\begin{abstract}
    We explore efficient estimation of statistical quantities, particularly rare event probabilities, for stochastic reaction networks. Consequently, we propose an importance sampling (IS) approach to improve the Monte Carlo (MC) estimator efficiency based on an approximate tau-leap  scheme. The crucial step in the IS framework is choosing an appropriate change of probability measure to achieve substantial variance reduction. This task is typically challenging and often requires insights into the underlying problem. Therefore, we propose an automated approach  to obtain a highly efficient path-dependent  measure change based on an original connection in the stochastic reaction network context between finding optimal IS parameters within a class of probability measures and a stochastic optimal control formulation. Optimal IS parameters are obtained by solving a variance minimization problem. First, we derive an associated dynamic programming equation. Analytically solving this backward equation is challenging, hence we propose an approximate dynamic programming formulation to find near-optimal control parameters. To mitigate the curse of dimensionality, we  propose a learning-based  method to approximate the value function using a neural network, where the parameters are determined via  a stochastic optimization algorithm. Our analysis and numerical experiments verify that the proposed learning-based  IS approach substantially reduces MC estimator variance, resulting in a lower computational complexity in the rare event regime, compared with standard tau-leap MC estimators.
    
    
\end{abstract}

\thispagestyle{plain}

\setcounter{tocdepth}{1}

\section{Introduction}
We propose an approach to efficiently estimate statistical quantities, particularly rare event probabilities for a particular class of continuous-time Markov chains known as stochastic reaction networks (SRNs). Consequently, we develop a learning-based importance sampling (IS) algorithm to improve the Monte Carlo (MC) estimator efficiency based on an approximate tau-leap (TL) scheme.  The automated approach is based on an original connection between optimal IS parameter determination within a class of probability measures  and stochastic optimal control (SOC) formulation.

SRNs (see Section \ref{sec:pjp} for a short introduction and \cite{ben2020hierarchical} for more details) describe the time evolution of biochemical reactions, epidemic processes \cite{brauer2001mathematical,anderson2015stochastic}, and transcription and translation in genomics and virus kinetics \cite{srivastava2002stochastic,hensel2009stochastic},  among other important applications. For the current study, let $\mathbf{X}$ be an SRN that takes values in $\nset^d$  and is defined in the time interval $[0,T]$, where $T>0$ is a user-selected final time. We aim to provide accurate and computationally efficient  MC estimations for the expected value $\mathbb{E}[g(\mathbf{X}(T))]$,  where $g:\nset^d\to\rset$ is a scalar observable for $\mathbf{X}$. In particular, we study estimating rare event probabilities with $g(\mathbf{x})=\bold{1}_{\{\mathbf{x} \in \mathcal{B}\}}$ (\ie, the indicator function for a set $\mathcal{B} \subset \rset^d$).

The quantity of interest, $\mathbb{E}[g(\mathbf{X}(T))]$, can be computed by solving the corresponding Kolmogorov backward equations~\cite{bayer2016efficient}. For most SRNs, deriving a closed-form solution for these ordinary differential equations is infeasible, and numerical approximations based on discretized schemes are commonly used. However, the computational cost scales exponentially with the number of species $d$. Therefore, we are  particularly interested in estimating $\mathbb{E}[g(\mathbf{X}(T))]$ using MC methods, an attractive alternative to avoid the curse of dimensionality. 

Many schemes have been developed to simulate exact sample paths for SRNs, such as the {stochastic simulation algorithm}~\cite{gillespie1976general} and {modified next reaction method}~\cite{anderson2007modified}. Pathwise exact SRN realizations can incur high computational costs if any reaction channels have high reaction rates. Gillespie~\cite{gillespie2001approximate} and Aparicio and Solari~\cite{aparicio2001population} independently proposed the explicit TL method (see Section~\ref{sec:exp_tau})  to overcome this issue by simulating approximate paths of $\mathbf{X}$, evolving the process with fixed time steps and keeping reaction rates fixed within each time step. Various simulation schemes have been subsequently proposed to deal with situations incorporating well-separated fast and slow time scales~\cite{cao2005trapezoidal,rathinam2007reversible,abdulle2010chebyshev,ahn2013implicit,moraes2016multilevel_splitting,hammouda2017multilevel}. 

Various variance reduction techniques have been proposed in the SRN context to reduce the computational work to estimate $\mathbb{E}[g(\mathbf{X}(T))]$. Several multilevel Monte Carlo (MLMC)~\cite{giles2008multilevel,giles2015multilevel} based methods have been proposed to address specific challenges in this context~\cite{Anderson2012,lester2015adaptive,moraes2016multilevel,moraes2016multilevel_splitting,hammouda2017multilevel,hammouda2020importance}. Furthermore, as naive MC and MLMC estimators fail to efficiently and accurately estimate rare event probabilities, different IS approaches \cite{kuwahara2008efficient,gillespie2009refining,roh2010state,daigle2011automated,cao2013adaptively,gillespie2019guided,roh2019data} have been proposed. 

The current paper proposes a path-dependent IS approach based on an approximate TL scheme to improve the MC estimator efficiency, and hence efficiently estimate various statistical quantities for SRNs (particularly rare event probabilities).  Our class of probability measure change is based on modifying the Poisson random variable rates used to construct the TL paths. In particular, optimal IS parameters are obtained by minimizing the second moment of the IS estimator (equivalently the variance) which represents the cost function for the associated SOC problem. We show that the corresponding value function solves  a dynamic programming relation that is challenging to solve analytically (see Section~\ref{sec:Approach Formulation}). We approximate the dynamic programming equation to derive a closed form solution and near-optimal control parameters. The cost to solve the associated backward equation numerically in multi-dimensional settings increases exponentially with respect to the dimension (\ie, the curse of dimensionality).  Thus, we propose approximating the resulting value function using a neural network to overcome this issue. Utilizing  the optimality criterion for the SOC problem, we obtain a relationship between  optimal IS parameters and the value function. Finally, we employ a stochastic optimization algorithm to learn the corresponding neural network parameters. Our analysis and numerical results for different dimensions confirm that the proposed estimator considerably reduces the variance compared with the standard TL-MC method with a negligible additional cost. This allows rare event probabilities to be efficiently computed in a regime where standard TL-MC estimators commonly fail.

The proposed approach is more computationally efficient than previously proposed IS schemes in this context~(\cite{kuwahara2008efficient,gillespie2009refining,roh2010state,daigle2011automated,cao2013adaptively,gillespie2019guided,roh2019data}) because it is based on an approximate TL scheme rather than the exact scheme. In contrast to previous approaches, the change  of measure is systematically  derived to ensure convergence to the optimal measure within the chosen class of probability measures, minimizing MC estimator variance. The novelty of this work is establishing a connection between IS and SOC in the context of pure jump processes, particularly for SRNs, with an emphasis on related practical and numerical  aspects. Note that some previous studies ~\cite{fleming2006controlled,dupuis2012importance,banisch2013meshless,hartmann2014characterization,zhang2014applications,hartmann2017variational,kebiri2017adaptive,hartmann2018importance,hartmann2019variational,nusken2021solving} have established a similar connection, mainly in the diffusion dynamics context, with less focus on pure jump dynamics. In this work, the proposed methodology is based on an approximate explicit TL scheme, which could and be subsequently extended in future work to continuous-time formulation (exact schemes),  and  implicit TL schemes which are relevant for systems with fast and slow time scales.

The remainder of this paper is organized as follows. Sections~\ref{sec:pjp},~\ref{sec:exp_tau},~\ref{sec:Monte Carlo (MC) estimator} and~\ref{sec:Importance Sampling} define relevant SRN, TL, MC and IS concepts, respectively. Section~\ref{sec:IS_SOC} establishes the connection between IS and SOC, formulating the SOC problem and defining its main ingredients: controls, cost function, and value function; then presents the dynamic programming solved by the optimal controls. Section~\ref{sec:Alternative Approach  for High Dimension} develops the proposed IS learning-based approach appropriate for multi-dimensional SRNs. Section \ref{sec:num_experiments} provides selected numerical experiments for different dimensions to illustrate the proposed approach's efficiency compared with standard MC approaches. Finally, Section~\ref{sec:Conclusions and future work} summarizes and concludes the work, and discusses possible future research directions.

\subsection{Stochastic reaction networks (SRNs)}
\label{sec:pjp}
We are interested in the time evolution for an homogeneously mixed  chemical reacting system described by the Markovian pure jump process,  $\mathbf{X}:[0,T]\times \Omega \to \nset^d$, where ($\Omega$, $\mathcal{F}$, $\mathbb{P}$)  is a probability space. In this framework, we assume that $d$ different species interact through $J$ reaction channels. 
The $i$-th component,  $X_i(t)$, describes the abundance of the $i$-th species present in the chemical system at time $t$. This work studies the time evolution of the state vector, 
\begin{equation}  
	\mathbf{X}(t) = \left(X_1(t), \ldots, X_d(t)\right) \in
	\nset^d \PERIOD
\end{equation}
Each reaction channel $\mathcal{R}_j$ is a pair $(a_j, \boldsymbol{\nu}_{j})$ defined by its propensity function $a_{j}:\rset^{d} \rightarrow \rset_{+}$ and stoichiometric vector $ \boldsymbol{\nu}_{j}=( \nu_{j,1},\nu_{j,2},..., \nu_{j,d})^\top$ satisfying
\begin{align}\label{reaction_channel}
	\mathbb{P}\left(\mathbf{X}(t+ \Delta t)=\mathbf{x}+ \boldsymbol{\nu}_{j} \mid  \mathbf{X}(t)=\mathbf{x}\right)=a_{j}(\mathbf{x})\Delta t + \ordo{\Delta t}, \: j=1,2,...,J  \PERIOD 
\end{align}
Thus, the probability of observing a jump in the process $\mathbf{X}$  from  state $\mathbf{x}$ to  state $\mathbf{x} +  \boldsymbol{\nu}_{j}$, a consequence of  reaction $\mathcal{R}_{j}$ firing during the small time interval $(t, t + \Delta  t]$, is proportional to the time interval length, $\Delta  t$, where $a_{j}(\mathbf{x})$  is the proportionality constant. We set $a_j(\mathbf{x}){=}0$ for  $\mathbf{x}$ such that $\mathbf{x}{+}\boldsymbol{\nu}_j\notin \nset^d$ (\ie, the \emph{non-negativity assumption}: the system can never produce negative population values). 

Hence, from \eqref{reaction_channel}, process $\mathbf{X}$  is a continuous-time, discrete-space Markov chain that can be characterized by Kurtz's random time change representation \cite{kurtz_2005},
\begin{equation}
	\label{eq:exact_process}
	\mathbf{X}(t)= \mathbf{x}_{0}+\sum_{j=1}^{J} Y_j \LP \int_0^t  a_{j}(\mathbf{X}(s)) \, \ud s \RP \boldsymbol{\nu}_j  \COMMA
\end{equation}
where $Y_j:\rset_+{\times} \Omega \to \nset$ are independent unit-rate Poisson processes. 
Conditions on the  reaction channels  can be imposed to ensure uniqueness~\cite{anderson2015stochastic} and avoid explosions in finite time~\cite{engblom2012stability,rathinam2013moment,gupta2014scalable}.

Applying the \textit{stochastic mass-action kinetics} principle,  we can assume that the propensity function $a_j(\cdot)$ for reaction channel $\mathcal{R}_j$, represented as\footnote{$\alpha_{j,i}$ molecules for species $S_i$ are consumed and $\beta_{j,i}$ are produced. Thus, $(\alpha_{j,i},\beta_{j,i}) \in \nset^2$ but $\beta_{j,i}-\alpha_{j,i}$ can be a negative integer, constituting the vector $\boldsymbol{\nu}_j=\left(\beta_{j,1}-\alpha_{j,1},\dots,\beta_{j,d}-\alpha_{j,d}\right) \in \zset^d$.}
\begin{equation}
	\alpha_{j,1} S_1+\dots+\alpha_{j,d} S_d \overset{\theta_j}{\rightarrow}\beta_{j,1} S_1+\dots+\beta_{j,d} S_d 
\end{equation}
obeys 
\begin{equation}\label{eq:prop_dynamics}
	a_j(\mathbf{x}):=\theta_j \prod_{i=1}^d \frac{x_i!}{(x_i-\alpha_{j,i})!} \mathbf{1}_{\{x_i\ge \alpha_{j,i}\}}\COMMA
\end{equation}
where $\{\theta_j\}_{j=1}^J$ represents positive constant reaction rates, and $x_i$ is the counting number for species $S_i$.

\subsection{Explicit tau-leap approximation}
\label{sec:exp_tau}
The explicit-TL scheme is a pathwise approximate method~\cite{gillespie2001approximate,aparicio2001population} to overcome computational drawbacks for exact methods (\ie, when many reactions fire during a short time interval). This scheme can be derived from the random time change representation~\eqref{eq:exact_process} by approximating the integral $\int_{t_i}^{t_{i+1}} a_{j}(\mathbf{X}(s)) \ud s $ as $a_j(\mathbf{X}(t_i))\,(t_{i+1}-t_i)$, \ie, using the forward-Euler method with time mesh  $\{t_{0}=0, t_{1},...,t_{N}= T\}$ and size $\Delta t=\frac{T}{N}$. Thus, the  \nameexp approximation for $\mathbf{X}$ should satisfy for $k\in\{1,2,\ldots,N\}$
\begin{equation}
	\hat{\mathbf{X}}^{\Delta t}_k = \mathbf{x}_{0}+\sum_{j=1}^{J} Y_{j} \LP  \sum_{i=0}^{k-1} a_{j}(\hat{\mathbf{X}}^{\Delta t}_i) \Delta t \RP   \boldsymbol{\nu}_{j} \COMMA
\end{equation}
and given $\hat{\mathbf{X}}_0:= \mathbf{x}_{0}$, we iteratively simulate a path for $\hat{\mathbf{X}}^{\Delta t}$ as 
\begin{equation}\label{eq:TL_approx}
	\hat{\mathbf{X}}^{\Delta t}_k:=\hat{\mathbf{X}}^{\Delta t}_{k-1}+\sum_{j=1}^{J} \mathcal{P}_{k-1,j}\left(a_{j}(\hat{\mathbf{X}}^{\Delta t}_{k-1}) \Delta t\right)  \boldsymbol{\nu}_{j} \COMMA \: 1 \le k \le N,
\end{equation}
where, conditioned on the current state $\hat{\mathbf{X}}^{\Delta t}_{k}$, $\{\mathcal{P}_{k,j}(r_{k,j})\}_{\{1\leq j\leq J \}}$ are independent Poisson random variables with respective rates $r_{k,j}:=a_{j}(\hat{\mathbf{X}}^{\Delta t}_{k})\Delta t$.

The  explicit-TL path $\hat{\mathbf{X}}^{\Delta t}$ is defined only at time mesh points, but can be naturally extended to $[0,T]$  as a piecewise constant path. We  apply the projection to zero to prevent  the process from exiting the lattice (\ie, producing negative values), hence \eqref{eq:TL_approx} becomes
\begin{equation}\label{eq:TL_approx_projection}
	\hat{\mathbf{X}}^{\Delta t}_k:=\max \left(\textbf{0},\hat{\mathbf{X}}^{\Delta t}_{k-1}+\sum_{j=1}^{J} \mathcal{P}_{k-1,j}\left(a_{j}(\hat{\mathbf{X}}^{\Delta t}_{k-1}) \Delta t\right)  \boldsymbol{\nu}_{j} \right) \COMMA \: 1 \le k \le N,
\end{equation}
where the maximum is applied entry-wise.
In this work, we use uniform time steps with length $\Delta t$, but the explicit-TL scheme and the proposed IS scheme (see Section \ref{sec:IS_SOC}) can also be applied to non-uniform time meshes.

\subsection{Biased Monte Carlo estimator}  
\label{sec:Monte Carlo (MC) estimator}
Let $\mathbf{X}$ be a stochastic process and $g: \rset ^{d} \rightarrow \rset$ a  scalar observable. 
We want to approximate $\mathbb{E} \left[g(\mathbf{X}(T))\right]$, but rather than sampling directly from $\mathbf{X}(T)$, we sample from $\overline{\mathbf{X}}^{\Delta t}(T)$, which are random variables generated by a numerical scheme  with step size $\Delta t$. 
We assume  that variates $\overline{\mathbf{X}}^{\Delta t}(T)$ are generated with an algorithm with weak order, $\Ordo{\Delta t}$, \ie, for sufficiently small $\Delta t$,
\begin{align}\label{eq:weak}
	\left |\mathbb{E} \left[g(\mathbf{X}(T))- g(\overline{\mathbf{X}}^{\Delta t}(T) )\right]\right |\leq C\Delta t
\end{align}
where $C>0$.\footnote{Refer to \cite{li2007analysis} for the underlying assumptions and proofs for this statement in the TL scheme context.}

Let $\mu_{M}$ be the standard MC estimator  for $\mathbb{E} \left[g(\overline{\mathbf{X}}^{\Delta t}(T))\right]$, 
\begin{equation}
	\mu_{M} :=\frac{1}{M}\sum_{m=1}^{M} g(\overline{\mathbf{X}}^{\Delta t}_{[m]}(T))\COMMA
\end{equation}
where $\{\overline{\mathbf{X}}^{\Delta t}_{[m]}(T)\}_{m=1}^M$ are independent and distributed as $\overline{\mathbf{X}}^{\Delta t}(T)$. 

The global error for the proposed MC estimator has error decomposition
\begin{align}\label{eq:error_split1}
	\left|\mathbb{E}[g(\mathbf{X}(T))]-\mu_M\right|\leq \underbrace{\left|\mathbb{E}[g(\mathbf{X}(T))]-\mathbb{E}[g(\overline{\mathbf{X}}^{\Delta t}(T))]\right|}_{\text{Bias}}+\underbrace{\left|\mathbb{E}[g(\overline{\mathbf{X}}^{\Delta t}(T))]-\mu_M\right|}_{\text{Statistical Error}}.
\end{align}
To achieve the desired accuracy, $\text{TOL}$, it is sufficient to bound the bias and statistical error equally by $\frac{TOL}{2}$. From \eqref{eq:weak}, choosing step size 
\begin{align}\label{eq:dtstar}
	\Delta t(\text{TOL})= \frac{\text{TOL}}{2\cdot C} 
\end{align}
ensures a bias of $\frac{\text{TOL}}{2}$.

Thus, considering the central limit theorem, the statistical error can be approximated as
\begin{align}\label{eq:staterror}
	|\mathbb{E}[g(\overline{\mathbf{X}}^{\Delta t}(T))]-\mu_{M}| \approx    C_{\alpha}\cdot \sqrt{\frac{\text{Var}[g(\overline{\mathbf{X}}^{\Delta t}(T))]}{M}},
\end{align}
where constant $C_{\alpha}$ is the $(1-\frac{\alpha}{2})-$quantile for the standard normal distribution. We choose $C_{\alpha}=1.96$ for a $95\%$ confidence level corresponding to $\alpha =0.05$.  Choosing 
\begin{align}\label{eq:Mstar1}
	M^*(\text{TOL})=C_{\alpha}^2\frac{4\cdot \text{Var}[g(\overline{\mathbf{X}}^{\Delta t}(T))]}{\text{TOL}^2}
\end{align} 
sample paths ensures the statistical error to be approximately bounded by $\frac{\text{TOL}}{2}$.

Given that the computational cost to simulate a single path is $\Ordo{ {\Delta t}^{-1}}$,  the expected total computational complexity is $\Ordo{\text{TOL}^{-3}}$; and the complexity scales with $ \text{Var}[g(\overline{\mathbf{X}}^{\Delta t}(T))]$ (see $\eqref{eq:Mstar1}$).

\subsection{Importance sampling}  
\label{sec:Importance Sampling}
Importance sampling (IS) techniques improve the computational costs for the crude MC estimator by variance reduction when used appropriately. To motivate the use of these techniques, consider estimating rare event probabilities, where the crude MC method is substantially expensive. In particular, consider estimating $q=\mathbb{P}(Y>\gamma)=\mathbb{E}[\textbf{1}_{\{Y>\gamma\}}]$, where $Y$ is a random variable taking values in $\mathbb{R}$ with probability density function  $\rho_{Y}$. Let $\gamma$ be sufficiently large that $q$ becomes sufficiently small. We can approximate $q$ using the MC estimator
\begin{align}
	\hat{q}=\frac{1}{M}\sum_{i=1}^M \textbf{1}_{\{Y^{(i)}>\gamma\}},
\end{align}
where $\{Y^{(i)}\}_{i=1}^{M}$ are independent and identically distributed (i.i.d)  realizations sampled according to $\rho_Y$. The MC estimator variance is 
\begin{align}\label{eq:varrareevent}
	Var\left[\mathbf1_{\{Y^{(i)}>\gamma\}}\right]&=q-q^2.
\end{align}
For a sufficiently small $q$, we can use \eqref{eq:varrareevent} and the central limit theorem to approximate the relative error as
\begin{align}
	\frac{|q-\hat{q}|}{q}\approx C_{\alpha}\sqrt{\frac{1}{qM}},
\end{align}
where $C_{\alpha}$ is chosen as in \eqref{eq:staterror}.

The number of required samples to attain a relative error tolerance $TOL_{rel}$ is $M\approx \frac{C_{\alpha}^2}{q\cdot TOL_{rel}^2}$. Thus, for $q$ of the order of $10^{-8}$, the number of required samples such that $TOL_{rel} = 5\%$  is approximately equal to $1.5\cdot 10^{11}$.

To demonstrate the IS concept, consider the general problem of estimating  $\mathbb{E}[g(Y)]$, where $g$ is a given observable. In the previous example, $g$ was chosen as $g(y)=\mathbf1_{\{y>\gamma\}}$. Let $\hat{\rho}_Z$ be the probability density function for a new real random variable $Z$, such that $g\cdot\rho_Y$ is dominated by  $\hat{\rho}_Z$, \ie,
\begin{align}\label{eq:ISsupport}
	\hat{\rho}_Z(x)=0 \implies g(x)\cdot\rho_Y(x)=0
\end{align}
for all $x\in \mathbb{R}$. This permits, the quantity of interest to be expressed as 
\begin{align}
	\mathbb{E}[g(Y)]=\int_{\mathbb{R}}g(x)\rho_Y(x)dx=\int_{\mathbb{R}}g(x)\underbrace{\frac{\rho_Y(x)}{\hat{\rho}_Z(x)}}_{L(x)}\cdot \hat{\rho}_Z(x) dx=\mathbb{E}[L(Z)\cdot g(Z)],
\end{align}
where $L(\cdot)$ is the likelihood ratio. Hence the expected value under the new measure remains unchanged, but the variance could be reduced due to a different second moment $\mathbb{E}\left[\left(g(Z)\cdot L(Z)\right)^2\right]$. 

The MC estimator under the IS measure is 
\begin{align}
	\mu_{M}^{IS}=\frac{1}{M} \sum_{j=1}^M L(Z_{[j]})\cdot g(Z_{[j]})=\frac{1}{M} \sum_{j=1}^M \frac{\rho_Y(Z_{[j]})}{\hat{\rho}_Z(Z_{[j]})}\cdot g(Z_{[j]}),
\end{align}
where $Z_{[j]}$ are i.i.d  samples from $\hat{\rho}_Z$ for $j=1,\dots,M$.

The main challenge when using IS is choosing a new probability measure that substantially reduces the variance compared with the original measure.
This step strongly depends on the structure of the problem under consideration. Further, the new measure should be obtained with negligible computational cost to ensure a computational efficient IS scheme. This is particularly challenging in the present problem, since we are considering path-dependent probability measures. In particular, the aim is to introduce a path-dependent change of probability measure that corresponds to changing the Poisson random variable rates used to construct the TL paths. Section~\ref{sec:Approach Formulation} shows how the optimal IS parameters can be obtained using a novel connection with SOC.

\section{Importance Sampling (IS) via Stochastic Optimal Control (SOC)}
\label{sec:IS_SOC}

\subsection{ Dynamic programming for the importance sample parameters}
\label{sec:Approach Formulation}
This section, establishes the connection between optimal IS measure determination within a class of probability measures, and SOC. 
Let $\mathbf{X}$ be a SRN as defined in Section~\ref{sec:pjp} and let $\hat{\mathbf{X}}^{\Delta t}$ denote its TL approximation as given by \eqref{eq:TL_approx_projection}. We aim to find a near-optimal IS measure to improve the MC estimator computational performance to estimate $\mathbb{E} \left[g(\mathbf{X}(T))\right]$. Since finding the optimal path-dependent change of measure within all measure classes presents a challenging problem, we limit ourselves to a parameterized class obtained via modifying the Poisson random variable rates of the TL paths. This class of measure change was previously used in \cite{hammouda2020importance}  to improve the MLMC estimator robustness and performance in this context; we focus on a single-level MC setting, and seek to automate the task to find a near-optimal IS measure within this class.

We introduce the change of measure resulting from changing the Poisson random variable rates in the TL scheme,
\begin{equation}\label{eq:measure_change}
	\bar{P}_{n,j}=\mathcal{P}_{n,j}\left(\delta_{n,j}^{\Delta t}(\overline{\mathbf{X}}^{\Delta t}_n)\Delta t\right),~~~ n=0,\dots, N-1, j=1,\dots,J ; 
\end{equation}
where $\delta_{n,j}^{\Delta t}(\mathbf{x})\in\mathcal{A}_{\mathbf{x},j}$ is the control parameter at time step $n$, under reaction $j$, and in state $\mathbf{x}\in\mathbb{N}^d$; and conditioned on $\overline{\mathbf{X}}^{\Delta t}_{n}$, $\mathcal{P}_{n,j}(r_{n,j})$ are  independent Poisson random variables with respective rates $r_{n,j}:=\delta_{n,j}^{\Delta t}(\overline{\mathbf{X}}^{\Delta t}_n)\Delta t$. 
The admissible set, 
\begin{align}\label{eq:addmissibleset}
	\mathcal{A}_{\mathbf{x},j}=\begin{cases}
		\{0\}&,\text{if }a_j(\mathbf{x})=0\\
		\{y\in\mathbb{R}: y>0\}&,\text{otherwise},
	\end{cases}
\end{align} 
is chosen such that (\ref{eq:ISsupport}) is fulfilled and to avoid infinite variance for the IS estimator. The control $\delta_{n,j}^{\Delta t}(\mathbf{x})\in\mathcal{A}_{\mathbf{x},j}$ depends deterministically on the current time step $n$, reaction channel $j$, and current state $\mathbf{x}=\overline{\mathbf{X}}^{\Delta t}_n$ for the TL-IS approximation in \eqref{eq:path_IS}.

Therefore, the  resulting scheme under the new measure is 
\begin{align}\label{eq:path_IS}
	\overline{\mathbf{X}}_{n+1}^{\Delta t}&=\max\left(\textbf{0},\overline{\mathbf{X}}_{n}^{\Delta t}+\sum_{j=1}^J\bar{P}_{n,j}\boldsymbol{\nu}_j\right) ,~~~ n=0,\dots,N-1,\\
	\overline{\mathbf{X}}_{0}^{\Delta t}&=\mathbf{x}_0,\nonumber   ; 
\end{align}
and the  likelihood ratio\footnote{We refer to \cite{hammouda2020importance} (Section 4.1) for the likelihood factor derivation of a similar IS scheme.} at step $n$ associated with the new IS measure is 
\begin{align}\label{eq:stepwiselh}
	L_n(\bar{\mathbf{P}}_n,\boldsymbol{\delta}_n^{\Delta t}(\overline{\mathbf{X}}^{\Delta t}_n))
	&=\prod_{j=1}^J\exp\left(-(a_j(\overline{\mathbf{X}}_{n}^{\Delta t})-\delta_{n,j}^{\Delta t}(\overline{\mathbf{X}}^{\Delta t}_n))\Delta t\right)\left(\frac{a_j(\overline{\mathbf{X}}_{n}^{\Delta t})}{\delta_{n,j}^{\Delta t}(\overline{\mathbf{X}}^{\Delta t}_n)}\right)^{\bar{P}_{n,j}}\nonumber\\
	&=\exp\left(-\left(\sum_{j=1}^J a_j(\overline{\mathbf{X}}_{n}^{\Delta t})-\delta_{n,j}^{\Delta t}(\overline{\mathbf{X}}^{\Delta t}_n)\right)\Delta t\right)  \cdot \prod_{j=1}^J\left(\frac{a_j(\overline{\mathbf{X}}_{n}^{\Delta t})}{\delta_{n,j}^{\Delta t}(\overline{\mathbf{X}}^{\Delta t}_n)}\right)^{\bar{P}_{n,j}}  ; 
\end{align} 
where $\boldsymbol{\delta}_n^{\Delta t}(\mathbf{x}) \in \times_{j=1}^J \mathcal{A}_{\mathbf{x},j}$ are the IS parameters with $\left(\boldsymbol{\delta}_n^{\Delta t}(\mathbf{x})\right)_j=\delta_{n,j}^{\Delta t}(\mathbf{x}) $ and the Poisson realizations are denoted by $\bar{\mathbf{P}}_n$ with $\left(\bar{\mathbf{P}}_n\right)_j:=\bar{P}_{n,j}$ for $j=1,\dots,J$. 
Equation (\ref{eq:stepwiselh}) uses the convention that  $\frac{a_j(\overline{\mathbf{X}}_{n}^{\Delta t})}{\delta_{n,j}^{\Delta t}(\overline{\mathbf{X}}^{\Delta t}_n)}=1$, whenever $a_j(\overline{\mathbf{X}}_{n}^{\Delta t})=0$ and $\delta_{n,j}^{\Delta t}(\overline{\mathbf{X}}^{\Delta t}_n)=0$. From (\ref{eq:addmissibleset}), this results in a factor of one in the likelihood ratio for reactions with  $a_j(\overline{\mathbf{X}}_{n}^{\Delta t})=0$.

Therefore, the likelihood ratio for $\{\overline{\mathbf{X}}^{\Delta t}_n: n=0,\dots,N\}$ across one path is 
\begin{equation}\label{eq:likelihood}
	L\left(\left(\bar{\mathbf{P}}_0,\dots,\bar{\mathbf{P}}_{N-1}\right),\left(\boldsymbol{\delta}_0^{\Delta t}(\overline{\mathbf{X}}^{\Delta t}_0),\dots,\boldsymbol{\delta}_{N-1}^{\Delta t}(\overline{\mathbf{X}}^{\Delta t}_{N-1})\right)\right)=\prod_{n=0}^{N-1} L_n(\bar{\mathbf{P}}_n,\boldsymbol{\delta}_n^{\Delta t}(\overline{\mathbf{X}}^{\Delta t}_n)).
\end{equation}
This likelihood ratio completes the characterization for the proposed IS approach, and allows the quantity of  interest  with respect to  the new measure to be expressed as
\begin{equation}\label{eq:expis}
	\mathbb{E}[g(\hat{\mathbf{X}}^{\Delta t}_N)]=\mathbb{E}\left[L\left(\left(\bar{\mathbf{P}}_0,\dots,\bar{\mathbf{P}}_{N-1}\right),\left(\boldsymbol{\delta}_0^{\Delta t}(\overline{\mathbf{X}}^{\Delta t}_0),\dots,\boldsymbol{\delta}_{N-1}^{\Delta t}(\overline{\mathbf{X}}^{\Delta t}_{N-1})\right)\right)\cdot g(\overline{\mathbf{X}}^{\Delta t}_N)\right],
\end{equation}
with the expectation in the right-hand side  of \eqref{eq:expis} taken with respect to the dynamics   in  \eqref{eq:path_IS}. 

Hereinafter, we aim to determine optimal parameters $\{\boldsymbol{\delta}_n^{\Delta t}(\mathbf{x})\}_{n=0,\dots,N-1; \mathbf{x}\in\mathbb{N}^d}$ that minimize the second moment (and hence the variance) for the IS estimator, given that $\overline{\mathbf{X}}_0^{\Delta t}=\mathbf{x}_0$. To that end, we derive an associated SOC formulation. First we introduce the cost function for the proposed SOC  problem in Definition~ \ref{def:secondmoment}, then derive a dynamic programming equation in  Theorem~\ref{theo:exact_optival} that is satisfied by the \textit{value function} $u_{\Delta t}(\cdot,\cdot)$  in Definition \ref{def:optival2}. The proof for Theorem \ref{theo:exact_optival} is given in Appendix~\ref{appendix:Proof of Theorem dynamic prog}.

\begin{definition}[Second moment for the proposed importance sampling estimator]\label{def:secondmoment}
	Let $0 \le n \le N$. Given that $\overline{\mathbf{X}}_{n}^{\Delta t}= \mathbf{x}$, the second moment for the proposed IS estimator can be expressed as
	\begin{align}\label{eq:second_moment}
		C_{n,\mathbf{x}}\left(\boldsymbol{\delta}^{\Delta t}_n,\dots,\boldsymbol{\delta}^{\Delta t}_{N-1}\right)&=\mathbb{E}\left[g^2(\overline{\mathbf{X}}_N^{\Delta t})\prod_{k=n}^{N-1} L_k^2 \left(\bar{\mathbf{P}}_k,\boldsymbol{\delta}_k^{\Delta t}(\overline{\mathbf{X}}_k^{\Delta t})\right) \middle| \overline{\mathbf{X}}_n^{\Delta t}=\mathbf{x}\right], \:  0 \le n \le N-1,
	\end{align}
	with terminal cost $C_{N,\mathbf{x}}=\mathbb{E}\left[g^2\left(\overline{\mathbf{X}}_N^{\Delta t}\right) \middle| \overline{\mathbf{X}}_N^{\Delta t}=\mathbf{x}\right]=g^2(\mathbf{x})$, for any $\mathbf{x} \in \mathbb{N}^d$. 
\end{definition}
Compared with the classical SOC formulation, \eqref{eq:second_moment} can be interpreted as the expected total cost; where the main difference is that \eqref{eq:second_moment} uses a multiplicative cost structure rather than the standard additive one. Therefore, we derive a dynamic programming relation in Theorem~\ref{theo:exact_optival} associated with this cost structure that is fulfilled by the corresponding value function (see Definition \ref{def:optival2}),  in the SRN context.

\begin{remark}[Structure of the cost function]
	One can derive an optimal control formulation with additive structure (similar to \cite{hartmann2017variational} in the stochastic differential equation setting) by applying a logarithmic transformation together with Jensen's inequality to \eqref{eq:second_moment}. 
	This reduces the control problem to a Kullback-Leibler minimization. In \cite{nusken2021solving,rached2022double}, this Kullback-Leibler minimization problem leads to the same optimal change of measure as the problem of finding the change of measure using a variance minimization approach.
	However, the previous conclusion needs more investigation in the setting of SRNs, which we leave for future potential work.
\end{remark}
\begin{definition}[Value function]\label{def:optival2}
	The \textit{value function} $u_{\Delta t}(\cdot,\cdot)$ is defined as  the optimal (infimum) second moment for the proposed IS estimator. For time step $0 \le n \le N$ and state $\mathbf{x} \in \mathbb{N}^d$, 
	\begin{align}
		u_{\Delta t}(n,\mathbf{x})&:=\inf_{\{\boldsymbol{\delta}^{\Delta t}_k\}_{k=n,\dots,N-1} \in \mathcal{A}^{N-n}}C_{n,\mathbf{x}}\left(\boldsymbol{\delta}^{\Delta t}_n,\dots,\boldsymbol{\delta}^{\Delta t}_{N-1}\right)\nonumber\\
		&=\inf_{\{\boldsymbol{\delta}^{\Delta t}_k\}_{k=n,\dots,N-1} \in \mathcal{A}^{N-n}}\mathbb{E}\left[g^2\left(\overline{\mathbf{X}}_N^{\Delta t}\right)\prod_{k=n}^{N-1} L_k^2\left(\bar{\mathbf{P}}_k,\boldsymbol{\delta}_k^{\Delta t}(\overline{\mathbf{X}}_k^{\Delta t})\right)\middle| \overline{\mathbf{X}}_n^{\Delta t}=\mathbf{x}\right],
	\end{align}
	where $\mathcal{A}=\bigtimes_{\mathbf{x}\in\mathbb{N}^d}\bigtimes_{j=1}^J\mathcal{A}_{\mathbf{x},j}\in\mathbb{R}^{\mathbb{N}^d \times J}$ is the admissible set for the IS parameters; and $u_{\Delta t}(N,\mathbf{x})=g^2(\mathbf{x})$, for any $\mathbf{x} \in  \mathbb{N}^d$.
\end{definition}
\begin{theorem}[Dynamic programming for importance sampling parameters]\label{theo:exact_optival}
	For $\mathbf{x}\in \mathbb{N}^d$, the value function $u_{\Delta t}(n,\mathbf{x})$ fulfills the dynamic programming relation
	\begin{small}
		\begin{align}\label{eq:exact_optival}
			u_{\Delta t}(N,\mathbf{x})&=g^2(\mathbf{x})\nonumber\\
			\text{and for } n&=N-1,\dots,0, \:  \text{and} \:  \mathcal{A}_\mathbf{x}:=\bigtimes_{j=1}^J\mathcal{A}_{\mathbf{x},j},\nonumber\\
			u_{\Delta t}(n,\mathbf{x})&=\inf_{\boldsymbol{\delta}_n^{\Delta t}(\mathbf{x})\in\mathcal{A}_\mathbf{x}}\exp\left(\left(-2\sum_{j=1}^J a_j(\mathbf{x})+\sum_{j=1}^J\delta_{n,j}^{\Delta t}(\mathbf{x})\right)\Delta t\right) \\
			&~~~~~~~~~~~~~\times \sum_{\mathbf{p} \in \mathbb{N}^J}\left(\prod_{j=1}^{J} \frac{(\Delta t \cdot \delta_{n,j}^{\Delta t}(\mathbf{x}))^{p_j}}{p_j!} (\frac{a_j(\mathbf{x})}{\delta_{n,j}^{\Delta t}(\mathbf{x})})^{2p_j} \right)\cdot u_{\Delta t}(n+1,\max(\mathbf{0},\mathbf{x}+ \boldsymbol{\nu}\mathbf{p})),\nonumber
		\end{align}
	\end{small}
	where $\boldsymbol{\nu}=\left(\boldsymbol{\nu}_1, \dots,\boldsymbol{\nu}_J\right)\in\mathbb{Z}^{d\times J}$.
\end{theorem}
Theorem~\ref{theo:exact_optival} breaks down the minimization problem to a simpler optimization that can be solved stepwise backward in time starting from final time $T$. Solving the minimization problem  \eqref{eq:exact_optival} analytically is difficult due to the infinite sum. Section \ref{sec:Algorithm} shows how to overcome this issue by approximating \eqref{eq:exact_optival} to derive near-optimal  parameters for $\{\boldsymbol{\delta}_n^{\Delta t}(\mathbf{x})\}_{n=0,\dots,N-1; \mathbf{x}\in\mathbb{N}^d}$ for the proposed IS approach.

\subsection{Approximate dynamic programming }
\label{sec:Algorithm}
Theorem \ref{theo:exact_optival} gives an exact solution for optimal IS parameters resulting from modifying the Poisson random variable rates in the TL paths. However, the infinite sum has to be evaluated in closed form to solve (\ref{eq:exact_optival}) analytically, which is generally difficult. Therefore, we propose approximating the value function $u_{\Delta t}(n,\mathbf{x})$ in (\ref{eq:exact_optival}) by $\overline{u}_{\Delta t}(n,\mathbf{x})$ for all time steps $n=0,\dots,N$, reaction channels $j=1,\dots,J$ and states $\mathbf{x}\in \mathbb{N}^d$. First, both $u_{\Delta t}(n,\mathbf{x})$ and $\overline{u}_{\Delta t}(n,\mathbf{x})$ satisfy the same final condition,
\begin{align}
	\overline{u}_{\Delta t}(N,\mathbf{x})=u_{\Delta t}(N,\mathbf{x})&=g^2(\mathbf{x}).
\end{align}

Next, to derive the approximate dynamic programming relation for $\overline{u}_{\Delta t}(\cdot,\cdot)$,  we presume Assumption \ref{assump:SOC_finite} to hold. This assumption is motivated by the behavior of the original propensities, which are of $\Ordo{1}$ due to the  mass-action kinetics principle (refer to \eqref{eq:prop_dynamics}).
\begin{assumption}\label{assump:SOC_finite}
	The controls $\{\boldsymbol{\delta}_{n}^{\Delta t}\}_{n=0,\dots,N-1}$ are asymptotically constant (\ie, $\delta_{n,j}^{\Delta t}(\mathbf{x}) \rightarrow c_{n,j,\mathbf{x}}$, as $\Delta t \rightarrow 0$, where $c_{n,j,\mathbf{x}}$ are constants for $1\leq j \leq J$, $0\leq n \leq N-1$, and $\mathbf{x}\in\mathbb{N}^d$).
\end{assumption}

Given Assumption  \ref{assump:SOC_finite} and that $\{a_j(\cdot)\}_{j=1}^J$ are of $\Ordo{1}$, we apply a Taylor expansion  around $\Delta t=0$  to the exponential term in \eqref{eq:exact_optival}, then truncate the expression within the infimum such that the remaining terms are $\Ordo{\Delta t}$. This truncates the infinite sum and linearizes the exponential term.
Thus, for $\mathbf{x}\in\mathbb{N}^d$ and $n=N-1,\dots,0$
\begin{small}
	\begin{align}\label{eq:approx_u}
		\overline{u}_{\Delta t}(n,\mathbf{x})
		&= \Delta t \inf_{(\delta_1,\dots,\delta_J)\in \mathcal{A}_\mathbf{x}}\big[ \sum_{j=1}^J \frac{a_j^2(\mathbf{x})}{\delta_j}  \overline{u}_{\Delta t}(n+1,\max(0, \mathbf{x}+\nu_j))\nonumber + \overline{u}_{\Delta t}(n+1,\mathbf{x}) \sum_{j=1}^J \delta_j\big]\nonumber\\
		&+\overline{u}_{\Delta t}(n+1,\mathbf{x})-2\Delta t\cdot \overline{u}_{\Delta t}(n+1,\mathbf{x}) \cdot \sum_{j=1}^J a_j(\mathbf{x})\nonumber\\
		&= \Delta t \cdot\sum_{j=1}^J\underbrace{ \inf_{\delta_j\in\mathcal{A}_{\mathbf{x},j}}\big[ \frac{a_j^2(\mathbf{x})}{\delta_j}\cdot \overline{u}_{\Delta t}(n+1,\max(0, \mathbf{x}+\nu_j))+ \delta_j\cdot \overline{u}_{\Delta t}(n+1,\mathbf{x})\big]}_{=:Q^{\Delta t}(n,j,\mathbf{x})}\nonumber\\
		&+\left(1-2\Delta t \sum_{j=1}^J a_j(\mathbf{x})\right)\overline{u}_{\Delta t}(n+1,\mathbf{x}),
	\end{align}
\end{small}
where $\delta_j \in \mathcal{A}_{\mathbf{x},j}$, $j=1,\dots,J$, are the SOC parameters at state $\mathbf{x}$ for reaction $j$. The admissible set $\mathcal{A}_{\mathbf{x},j}$ is defined in \eqref{eq:addmissibleset}.  Assumption \ref{assump:SOC_finite} ensures that (i) we can apply the Taylor expansion to the exponential term as $\Delta t$ decreases, and (ii) we have the exact approximation structure for \eqref{eq:approx_u} with no further terms scaling with $\Delta t$ that have order less than $\Delta t^2$.

The infimum in \eqref{eq:approx_u} is attained when
\begin{equation}\label{eq:condition of standard case}
	(i) \: \overline{u}_{\Delta t}(n+1,\mathbf{x}) \neq 0, \quad \text{and} \quad  (ii) \:  \overline{u}_{\Delta t}(n+1, \max(0,\mathbf{x}+\nu_j)) \neq 0, \: \forall \:  1 \le j \le J.
\end{equation}
In this case,  the approximate optimal SOC parameter $\overline{\delta}^{\Delta t}_{n,j}(\mathbf{x})$ can be analytically determined as 
\begin{align}\label{eq: standard_control_solution}
	\overline{\delta}^{\Delta t}_{n,j}(\mathbf{x})&= \frac{a_j(x)\sqrt{\overline{u}_{\Delta t}(n+1, \max(0,\mathbf{x}+\nu_j))}}{\sqrt{\overline{u}_{\Delta t}(n+1,\mathbf{x})}},\: 1 \le j \le J.
\end{align}

Note \eqref{eq: standard_control_solution} includes the particular case when  $a_j(\mathbf{x})=0$ for some $j\in\{1,\dots,J\}$. In such a case, $\overline{\delta}_{n,j}^{\Delta t}(\mathbf{x})=0$, which agrees with \eqref{eq:addmissibleset}.

An important advantage for this numerical approximation, $\overline{u}_{\Delta t}(\cdot,\cdot)$, is that we reduce the complexity of the original optimization problem at each step in \eqref{eq:exact_optival} from a  simultaneous optimization over $J$ variables to independent one-dimensional optimization problems that can be solved in parallel using \eqref{eq: standard_control_solution}. 

\begin{remark}[Assumption \eqref{eq:condition of standard case}]
	Whether the assumption in \eqref{eq:condition of standard case} is generally fulfilled depends on the method employed to solve the dynamic programming principle in \eqref{eq:approx_u}. For example, if we use a direct numerical implementation either some special numerical treatment is required for the cases where \eqref{eq:condition of standard case} is violated, or some regularization is required to ensure well-posedness. The proposed approach from Section~\ref{sec:Alternative Approach  for High Dimension} avoids that issue since we model $\overline{u}_{\Delta t}(\cdot,\cdot)$ with a strictly positive ansatz function, which guarantees  condition \eqref{eq:condition of standard case} to hold for any state $\mathbf{x}$ and all time steps $n$.
\end{remark}

\begin{remark}[Computational cost for dynamic programming]\label{rem:cursed}
	To derive a practical numerical algorithm for a finite number of states, we truncate the infinite state space $\mathbb{N}^d$ to $\bigtimes_{i=1}^d [0,\overline{S}_i]$, where $\overline{S}_1,\dots,\overline{S}_d$ is a set of sufficiently large upper bounds. 
	The computational cost to numerically solve the dynamic programming equation \eqref{eq:approx_u} for step size $\Delta t$ and state space $\bigtimes_{i=1}^d [0,\overline{S}_i]$ can be expressed as 
	\begin{align}\label{eq:costbackward}
		W_{\text{dp}}(\bar{\mathbf{S}},\Delta t)\approx  \left(\bar{S}^\ast\right)^d \cdot \frac{T}{\Delta t}\cdot J,
	\end{align}
	where  $\bar{S}^\ast=\max_{i=1,\dots,d}\bar{S}_i$.
	
	The cost in \eqref{eq:costbackward} scales exponentially with dimension $d$. Section~\ref{sec:Alternative Approach  for High Dimension} proposes an alternative approach to address this curse of dimensionality. However, in future work, we aim to combine dimension reduction techniques for SRNs with a direct numerical implementation of dynamic programming.
\end{remark}

\subsection{Learning-based  approach}
\label{sec:Alternative Approach  for High Dimension}

Using the SOC formulation derived in Section \ref{sec:Algorithm}, we propose approximating the value function  $\overline{u}_{\Delta t}(\cdot,\cdot)$ with a parameterized ansatz function, $\hat{u}(t,\mathbf{x};\boldsymbol{\beta})$.

\begin{remark}[Choosing the ansatz function]
	The parameterized ansatz function $\hat{u}(t,\mathbf{x};\boldsymbol{\beta})$ should consider the final condition of the value function \eqref{eq:exact_optival}, and its choice depends on the given SRN and observable $g(\mathbf{x})$. 
	For linear observables, such as  $g(\mathbf{x})=x_i$, we  can consider polynomial basis functions as an ansatz. For more complex problems, the ansatz function is a small neural network.
\end{remark}
For rare event applications with observable $g(\mathbf{\textbf{x}})=\mathbf{1}_{\{x_i>\gamma\}}$, we consider a sigmoid  with learning parameters $\boldsymbol{\beta}=\left(\boldsymbol{\beta}^{space},\beta^{time}\right) \in \mathbb{R}^{d+1}$ as the ansatz function
\begin{align}\label{eq:ansatz_sigmoid}
	\hat{u}(t,\mathbf{x};\boldsymbol{\beta})= \frac{1}{1+e^{-(1-t) \cdot \left(\langle\boldsymbol{\beta}^{space},\mathbf{x}\rangle+\beta^{time}\right)-b_0-\beta_0x_i}},
\end{align}
where $\langle\cdot,\cdot\rangle$ denotes  the inner product, and the time is scaled to one using $t\in[0,1]$.

Parameters  $b_0$ and $\beta_0$ are not learned through optimization but determined by  fitting  the final condition for Theorem \ref{theo:exact_optival}, which imposes $\hat{u}(1,\mathbf{x};\boldsymbol{\beta})\approx g^2(\mathbf{x})=\mathbf{1}_{\{x_i>\gamma\}}$. Therefore, the discontinuous indicator function is approximated by a sigmoid, and the fit is characterized by the position of the sigmoid's inflection point and the sharpness of the slope. The position and value of local and global minima with respect to the learned parameters $\boldsymbol{\beta}^{space}$ and ${\beta}^{time}$ depend on the choices for $b_0$ and ${\beta}_0$.

To derive  IS parameters from the ansatz function, we use the previous SOC result from \eqref{eq: standard_control_solution}, \ie,
\begin{align}\label{eq:deltafromu}
	\hat{\delta}^{\Delta t}_{j}(n,\mathbf{x};\boldsymbol{\beta})&= \frac{a_j(\mathbf{x})\sqrt{\hat{u}\left(\frac{(n+1)\Delta t}{T}, \max(0,\mathbf{x}+\nu_j);\boldsymbol{\beta}\right)}}{\sqrt{\hat{u}(\frac{(n+1)\Delta t}{T},\mathbf{x};\boldsymbol{\beta})}},\: 1 \le j \le J, \: 0 \le n\le N-1, \: \mathbf{x}\in \mathbb{N}^d.
\end{align}
We define $\hat{u}(t,\cdot;\cdot)$ in \eqref{eq:ansatz_sigmoid} as a time-continuous function for $t\in [0,1]$; whereas the IS controls from $\hat{\delta}^{\Delta t}_{j}(n,\cdot;\cdot)$ are discrete in time for $n=0,\dots,N-1$, and depend on time step size $\Delta t$. Therefore, $\hat{u}(\cdot,\cdot;\boldsymbol{\beta})$ can be used to derive control parameters for arbitrary $\Delta t$ in \eqref{eq:deltafromu}.

The parameters $\boldsymbol{\beta}$ for the ansatz function are then chosen to minimize the second moment,
\begin{align}\label{eq:optiansatz}
	\inf_{\boldsymbol{\beta\in \mathbb{R}^{d+1}}} \mathbb{E}\underbrace{\left[g^2\left(\overline{\mathbf{X}}_N^{\Delta t,\boldsymbol{\beta}}\right)\prod_{k=0}^{N-1} L_k^2\left(\bar{\mathbf{P}}_k,\hat{\boldsymbol{\delta}}^{\Delta t}(k,\overline{\mathbf{X}}_k^{\Delta t,\boldsymbol{\beta}};\boldsymbol{\beta})\right)\right]}_{=:C_{0,\mathbf{x}}\left(\hat{\boldsymbol{\delta}}^{\Delta t}_0,\dots,\hat{\boldsymbol{\delta}}^{\Delta t}_{N-1}; \boldsymbol{\beta}\right)},
\end{align}
where $\{\overline{\mathbf{X}}_n^{\Delta t,\boldsymbol{\beta}}\}_{n=1,\dots,N}$ is the IS path generated using IS parameters from \eqref{eq:deltafromu} and $\left(\hat{\boldsymbol{\delta}}^{\Delta t}(n,\mathbf{x};\boldsymbol{\beta})\right)_j=\hat{\delta}^{\Delta t}_{j}(n,\mathbf{x};\boldsymbol{\beta})$ for $1\leq j\leq J$.

We use a gradient based stochastic optimizer method to solve \eqref{eq:optiansatz}, and derive Lemma \ref{lem:gradient} (proof in Appendix~\ref{apdx:proofgradient}) for the gradient of the second moment with respect to parameters $\boldsymbol{\beta}$.

\begin{lemma}\label{lem:gradient}
	The partial derivatives for the second moment $C_{0,\mathbf{x}}\left(\hat{\boldsymbol{\delta}}^{\Delta t}_0,\dots,\hat{\boldsymbol{\delta}}^{\Delta t}_{N-1}; \boldsymbol{\beta}\right)$  in \eqref{eq:optiansatz} with respect to $\beta_{l}$, $l=1,\dots, (d+1)$, are  given by
	\begin{align}\label{eq:gradient}
		&\frac{\partial}{\partial \beta_l}\mathbb{E}\left[\underset{=:R(\mathbf{x}_0;\boldsymbol{\beta})}{\underbrace{g^2\left(\overline{\mathbf{X}}_N^{\Delta t,\boldsymbol{\beta}}\right)\prod_{k=0}^{N-1} L_k^2\left(\bar{\mathbf{P}}_k,\hat{\boldsymbol{\delta}}^{\Delta t}(k,\overline{\mathbf{X}}_k^{\Delta t,\boldsymbol{\beta}};\boldsymbol{\beta})\right)}}\right]\nonumber\\    
		&=\mathbb{E}\left[R(\mathbf{x}_0;\boldsymbol{\beta}) \left(\sum_{k=1}^{N-1}\sum_{j=1}^J \left(\Delta t - \frac{\bar{P}_{k,j}}{\hat{\delta}_j^{\Delta t}(k,\overline{\mathbf{X}}^{\Delta t,\boldsymbol{\beta}}_k;\boldsymbol{\beta})}\right)\cdot \frac{\partial}{\partial \beta_l} \hat{\delta}_j^{\Delta t}(k,\overline{\mathbf{X}}^{\Delta t,\boldsymbol{\beta}}_k;\boldsymbol{\beta})\right) \right],
	\end{align}
	where $\{\overline{\mathbf{X}}_n^{\Delta t,\boldsymbol{\beta}}\}_{n=1,\dots,N}$ is the IS path generated using the IS parameters from \eqref{eq:deltafromu} and
	\begin{small}
		\begin{align}\label{eq:gradientdelta}
			&   \frac{\partial}{\partial \beta_l} \hat{\delta}_j^{\Delta t}(k,\mathbf{x};\boldsymbol{\beta}) \\
			&=\frac{a_j^2(\mathbf{x})}{2\hat{\delta}_j^{\Delta t}(k,\mathbf{x};\boldsymbol{\beta})}\cdot\left(\frac{\frac{\partial }{\partial \beta_l}\hat{u}(\frac{(k+1)\Delta t}{T},\max(\mathbf{x}+\nu_j,0);\boldsymbol{\beta})}{\hat{u}(\frac{(k+1)\Delta t}{T},\mathbf{x};\boldsymbol{\beta})} -\frac{\hat{u}(\frac{(k+1)\Delta t}{T},\max(\mathbf{x}+\nu_j,0);\boldsymbol{\beta})\frac{\partial }{\partial \beta_l}\hat{u}(\frac{(k+1)\Delta t}{T},\mathbf{x};\boldsymbol{\beta})}{\hat{u}^2(\frac{(k+1)\Delta t}{T},\mathbf{x};\boldsymbol{\beta})}\right)\nonumber.
		\end{align}
	\end{small}
\end{lemma}

Thus, partial derivatives for $\hat{u}(t,\mathbf{x};\boldsymbol{\beta})$ for the ansatz \eqref{eq:ansatz_sigmoid} are 
\begin{small}
	\begin{align}\label{eq:partderivu}
		\frac{\partial }{\partial \beta_l}\hat{u}(t,\mathbf{x};\boldsymbol{\beta})
		&=\begin{cases}     (1-t)x_i\hat{u}(t,\mathbf{x};\boldsymbol{\beta})(1-\hat{u}(t,\mathbf{x};\boldsymbol{\beta}))&, \text{if }  \beta_l=\left(\boldsymbol{\beta}^{space}\right)_{i}\\
			(1-t)\hat{u}(t,\mathbf{x};\boldsymbol{\beta})(1-\hat{u}(t,\mathbf{x};\boldsymbol{\beta}))&, \text{if } \beta_l=\beta^{time},
		\end{cases}
	\end{align} 
\end{small} 
where $\left(\boldsymbol{\beta}^{space}\right)_{i}$ denotes the i-th entry for $\boldsymbol{\beta}^{space}$.

For an ansatz function different from \eqref{eq:ansatz_sigmoid}, the gradient is still given by Lemma \ref{lem:gradient} only the derivation of $\frac{\partial }{\partial \beta_l}\hat{u}(t,\mathbf{x};\boldsymbol{\beta})$ in \eqref{eq:partderivu} changes accordingly.

By estimating the gradient in  \eqref{eq:gradient} using a MC estimator, we iteratively optimize the parameters $\boldsymbol{\beta}$ to reduce the variance. For this optimization, we use the Adam optimizer with the same parameter values suggested in \cite{kingma2014adam} with the only difference that the step size is tuned to fit our problem setting.

In Section \ref{sec:num_experiments}, we illustrate the potential of our new IS method based on the learning approach numerically in terms of variance reduction. Further theoretical and numerical analysis of this approach is left for future work, particularly the initialization for the learned parameters $\beta^{time}$ and $\boldsymbol{\beta}^{space}$ in \eqref{eq:ansatz_sigmoid} and investigations of a stopping rule.

To derive an estimator for $\mathbb{E}[g(\mathbf{X}(T))]$ using the proposed IS change of measure, we first solve the related SOC problem using the approach from this section; then we simulate $M$ paths under the new IS sampling measure. Thus, the MC estimator using the proposed IS change of measure over $M$ paths becomes
\begin{align}\label{eq:ISMC}
	\mu^{IS}_{M,\Delta t}=\frac{1}{M} \sum_{i=1}^M L_i\cdot g(\overline{\mathbf{X}}_{[i],N}^{\Delta t,\boldsymbol{\beta}}),
\end{align}
where $\overline{\mathbf{X}}_{[i],N}^{\Delta t,\boldsymbol{\beta}}$ is the $i$-th IS sample path and the corresponding likelihood factor from \eqref{eq:likelihood} is 
\begin{align}\label{eq:lik_learning}
	L_i=L\left(\left(\bar{\mathbf{P}}_0,\dots,\bar{\mathbf{P}}_{N-1}\right),\left(\hat{\boldsymbol{\delta}}^{\Delta t}(0,\overline{\mathbf{X}}_{[i],0}^{\Delta t,\boldsymbol{\beta}};\boldsymbol{\beta}),\dots,\hat{\boldsymbol{\delta}}^{\Delta t}(N-1,\overline{\mathbf{X}}_{[i],N-1}^{\Delta t,\boldsymbol{\beta}};\boldsymbol{\beta})\right)\right) . 
\end{align} 

\begin{remark}
	The  explicit pathwise derivatives in  Lemma \ref{lem:gradient} have the following advantages compared with the finite difference approach: (i) the explicit pathwise derivatives  are unbiased with respect to the TL scheme,  resulting in only the MC error for evaluating the expectation (\ie, without additional finite difference error), and     
	(ii) the gradient computation in \eqref{eq:gradient} requires the estimation of an expected value with a high relative error because of $g$ being fitted to an indicator function. Using the IS-TL paths we control better the related statistical error.
\end{remark}

\subsection{Computational cost for the learning-based approach}
\label{sec:cost}
This section discusses the computational complexity for the learning approach to achieve a prescribed tolerance  $\text{TOL}$. Recall that the proposed approach comprises two steps; hence, two types of costs occur: (i) the offline learning cost for the ansatz function parameters $\boldsymbol{\beta}$, and 
(ii) the online cost to obtaining the MC estimator \eqref{eq:ISMC} based on $M$ simulated paths using the derived IS measure (see \eqref{eq:deltafromu}). 

The offline cost for (i) can be expressed as
\begin{align*}
	W_{pl}(I,M_0,\Delta t_{pl})\approx I \cdot M_0 \cdot \frac{T}{\Delta t_{pl}} \cdot J \cdot(C_{Poi} + C_{grad}),
\end{align*}
where $I$ is the number of optimizer steps, $M_0$ is the number of paths needed to derive the estimator of the gradient per optimizer step, $C_{Poi}$ is the cost to generate one Poisson random variable, $C_{grad}$ is the cost for the update of the algebraic evaluation of \eqref{eq:gradient}, and $\Delta t_{pl}$ is the step size. In contrast to \eqref{eq:costbackward}, this offline cost does not scale exponentially with dimension $d$.

The cost for one IS-TL path based on $\hat{u}(\cdot,\cdot;\boldsymbol{\beta})$ is the same as for a TL path with negligible additional factors $C_{\hat{\delta}}$ for evaluating \eqref{eq:deltafromu} and $C_{lik}$ for deriving the likelihood update, as given in \eqref{eq:stepwiselh},
\begin{align*}
	W_{forward}(\Delta t_f)\approx \frac{T}{\Delta t_f} \cdot J \cdot (C_{Poi} +C_{lik}+C_{\hat{\delta}}),
\end{align*}
where $\Delta t_f$ is the step size.
Thus, total cost is 
\begin{align*}
	W_{IS-TL}(M,\Delta t_{pl},\Delta t_f) \approx W_{pl} (I,M_0,\Delta t_{pl})+ M\cdot W_{forward}(\Delta t_f).
\end{align*}

Following the same derivation as for \eqref{eq:error_split1}-\eqref{eq:Mstar1}, we choose $\Delta t_{f}=\frac{\text{TOL}}{2 \cdot C}$, where $C$ is the constant from \eqref{eq:weak}, to obtain total computational complexity to derive a prescribed tolerance $TOL$ 
\begin{align}\label{eq:totalcost}
	W_{IS-TL}(\text{TOL})=W_{pl} (I,M_0,\Delta t_{pl})+const \cdot \frac{\text{Var}[g\left(\overline{\mathbf{X}}^{\Delta t,\boldsymbol{\beta}}_N\right)\cdot L]}{\text{TOL}^3},
\end{align}
where $L$ is the likelihood factor corresponding to the IS path $\overline{\mathbf{X}}^{\Delta t,\boldsymbol{\beta}}$ (refer to \eqref{eq:lik_learning}).

Our numerical simulations suggest that the amount of variance reduction achieved with the proposed approach is not related to $\Delta t_{pl}$ (see Figure~\ref{fig:4dII}). Therefore, we can achieve a low offline parameter learning cost ($W_{pl} (I,M_0,\Delta t_{pl})$) by using $\Delta t_{pl}\gg\Delta t_{f}$. 

For comparison, Section \ref{sec:Monte Carlo (MC) estimator} shows that the standard MC-TL approach has total computational complexity 
$$
W_{MC-TL}(TOL)=const_{TL}\cdot\frac{Var[g(\hat{\mathbf{X}}^{\Delta t}_N)]}{TOL^3} .
$$ 
The proposed IS approach reduces this cost by variance reduction $\left(\text{Var}[g(\overline{\mathbf{X}}^{\Delta t}_N)\cdot L]\ll Var[g(\hat{\mathbf{X}}^{\Delta t}_N)]\right)$ (refer to Figures~\ref{fig:decayadam}--\ref{fig:6d}). The TL variance becomes increasingly large in the asymptotic regime for very rare event probabilities, such that the additional cost $W_{pl} (I,M_0,\Delta t_{pl})$ for learning $\boldsymbol{\beta}$ in \eqref{eq:totalcost} becomes negligible. Therefore, we obtain $W_{IS-TL}(TOL)\ll W_{MC-TL}(TOL)$ in the rare event regime.

\section{Numerical Experiments and Results}
\label{sec:num_experiments}

Through Examples~\ref{exp:decay}, \ref{exp:mm}, and \ref{exp:6d}, we
demonstrate the advantages for the proposed IS approach compared with the standard MC approach. We numerically show that the proposed approach achieves substantial variance reduction compared with standard MC estimators when applied to SRNs with different dimensions.

\begin{example}[Pure decay]
	\label{exp:decay}
	This example considers one species and a single reaction,
	\begin{align*}
		X\overset{\theta_1}{\rightarrow} \emptyset,
	\end{align*}
	where $\theta_1=1$, and the final time $T = 1$. 
	Thus, the propensity is $a(x)=\theta_1x$, the stoichiometric vector is $\nu=-1$, and the observable is $g(x)=\mathbf1_{\{x>50\}}$ with $X_0=100$. 
\end{example}

\begin{example}[Michaelis-Menten enzyme kinetics] 
	\label{exp:mm}
	The Michaelis-Menten enzyme kinetics~\cite{rao2003stochastic} describe the catalytic conversion of substrate $S$ into a product $P$ through three reactions,
	\begin{align*}
		E+S\overset{\theta_1}{\rightarrow} C, ~~C\overset{\theta_2}{\rightarrow} E+S,~~
		C\overset{\theta_3}{\rightarrow} E+P,
	\end{align*}
	where $E$ denotes the enzyme and $\theta = (0.001,0.005,0.01)^\top$. 
	We consider the initial state $\mathbf{X}_0=(E(0),S(0),C(0),P(0))^\top=(100, 100, 0, 0)^\top$ and the final time $T=1$. The corresponding propensity and the change of the state matrix are 
	\begin{align*}
		a(\textbf{x})=\left(\begin{array}{c}
			\theta_{1} E S \\
			\theta_{2} C \\
			\theta_{3} C
		\end{array}\right), \quad
		\boldsymbol{\nu}=\left(\begin{array}{ccc}
			-1 & 1 & 1  \\
			-1 & 1 & 0  \\
			1 & -1& -1 \\
			0& 0& 1
		\end{array}\right).
	\end{align*}
	The observable of interest is $g(\mathbf{x})=\mathbf{1}_{\{x_3>22\}}$.
\end{example}

\begin{example}[Enzymatic futile cycle model]
	\label{exp:6d}
	The enzymatic futile cycle~\cite{kuwahara2008efficient} describes two instances for the elementary single-substrate enzymatic reaction scheme and can be described by six reactions,
	\begin{align*}
		&R_{1}: S_{1}+S_{2} \stackrel{\theta_{1}}{\longrightarrow} S_{3}\text {, } \quad R_{2}: S_{3} \stackrel{\theta_{2}}{\longrightarrow} S_{1}+S_{2} \text {, }\quad
		R_{3}: S_{3} \stackrel{\theta_{3}}{\longrightarrow} S_{1}+S_{5} \text {, }\\
		&R_{4}: S_{4}+S_{5} \stackrel{\theta_{4}}{\longrightarrow} S_{6} \text {, }\quad
		R_{5}: S_{6} \stackrel{\theta_{5}}{\longrightarrow} S_{4}+S_{5} \text {, }\quad
		R_{6}: S_{6} \stackrel{\theta_{6}}{\longrightarrow} S_{4}+S_{2} \text {. }
	\end{align*}
	Initial states are $\mathbf{X}(0)=\left(S_1(0),\dots,S_6(0)\right)=\left (1, 50, 0, 1, 50, 0 \right)$,  and we take the rates as $\theta_{1}=\theta_{2}=\theta_{4}=\theta_{5}=1$, and $\theta_{3}=\theta_{6}=0.1$. The propensity $a(\mathbf{x})$ follows the stochastic mass-action kinetics in \eqref{eq:prop_dynamics} and the final time is $T=2$. We consider $g(\mathbf{x})=\mathbf{1}_{\{x_5>60\}}$ as the observable.
\end{example} 

Since all three are rare event examples with observable $g(\mathbf{\textbf{x}})=\mathbf{1}_{\{x_i>\gamma\}}$, we use the ansatz function \eqref{eq:ansatz_sigmoid} with initial parameters $\boldsymbol{\beta}^{space}=0$, and $\beta^{time}=0$.
The relative error is more relevant for rare event occurrences than the absolute error, hence we use a relative version of the variance, \ie, the squared coefficient of variation~\cite{ben2021efficient, kroese2013handbook}, which, for a random variable $X$, is given by
\begin{align}
	Var_{rel}[X]=\frac{Var[X]}{\mathbb{E}[X]^2}.
\end{align}
To judge the robustness of our variance estimators, we estimate the kurtosis, $\kappa:=\frac{\mathrm{E}\left[\left(X-E\left[X\right]\right)^4\right]}{\left(\operatorname{Var}\left[X\right]\right)^2}$, because the standard deviation of the sample variance \cite{hammouda2020importance} is given by
\begin{align*}
	\sigma_{\mathcal{S}^2\left(X\right)}=\frac{\operatorname{Var}\left[X\right]}{\sqrt{M}} \sqrt{\left(\kappa-1\right)+\frac{2}{M-1}},
\end{align*}
where $M$ is the number of samples.

We set the Adam optimizer step size $\alpha=0.1$ for the simulations.

\begin{figure} [h]
	\caption{Example \ref{exp:decay} with step size $\Delta t_{pl}=\Delta t_f=1/2^4$ for the proposed IS-MC estimator: \textbf{(a)} sample mean;  \textbf{(b)} squared coefficient of variation; \textbf{(c)} parameters; \textbf{(d)} kurtosis for each optimizer step. Adam optimizer gradient, sample variance, and kurtosis were estimated using $M_0=10^4$ samples. The reference value for the standard MC-TL approach was derived from a single run with $M=10^6$ samples and with step size $\Delta t=1/2^4$.} 
	\label{fig:decayadam}
	\subfloat[]{\includegraphics[width=0.49\textwidth]{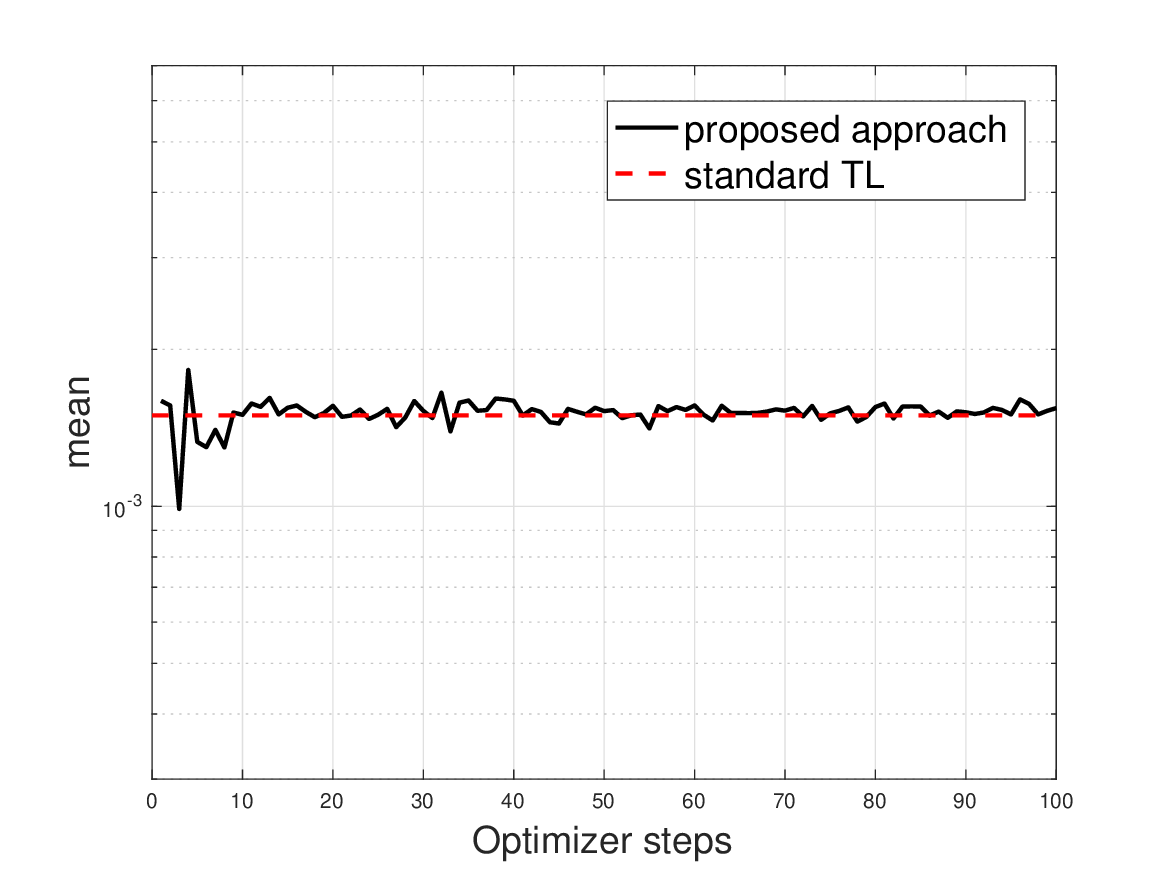}}
	\subfloat[]{\includegraphics[width=0.49\textwidth]{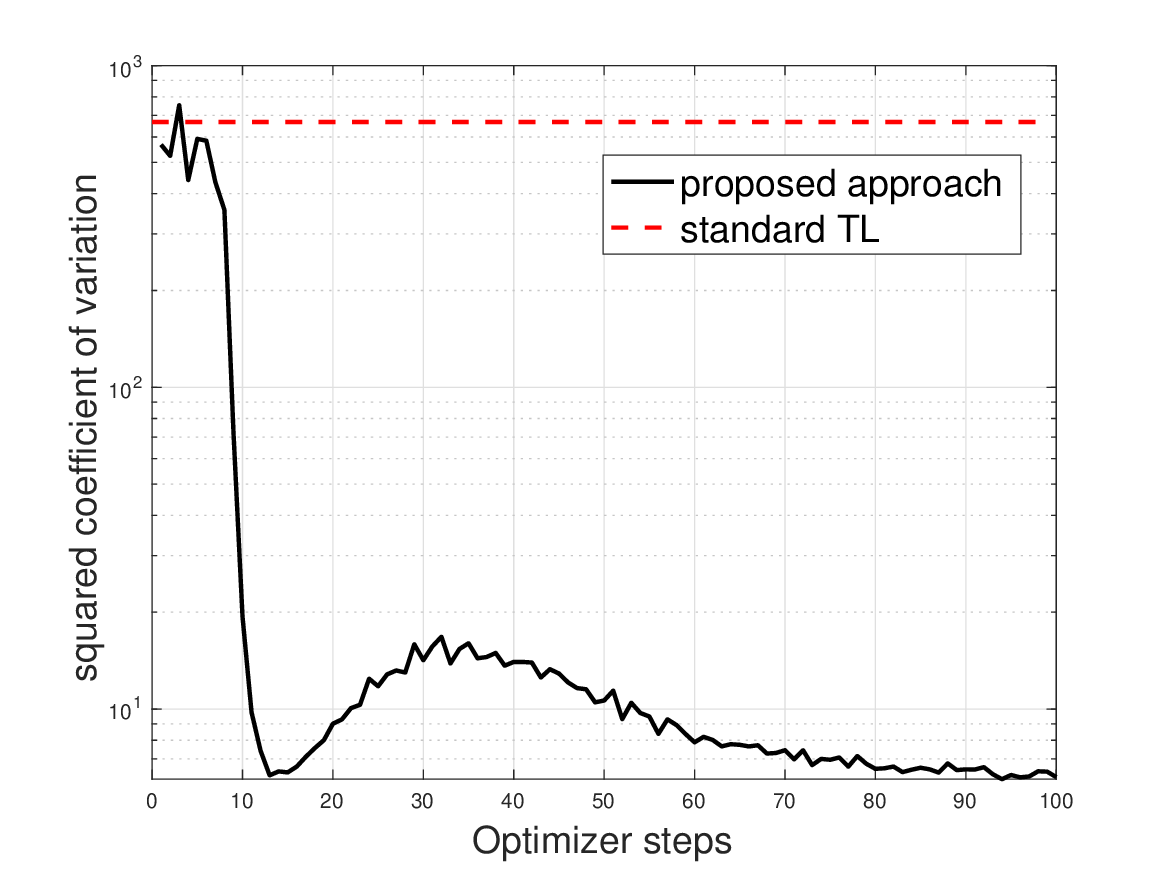}}  \\
	\subfloat[]{\includegraphics[width=0.49\textwidth]{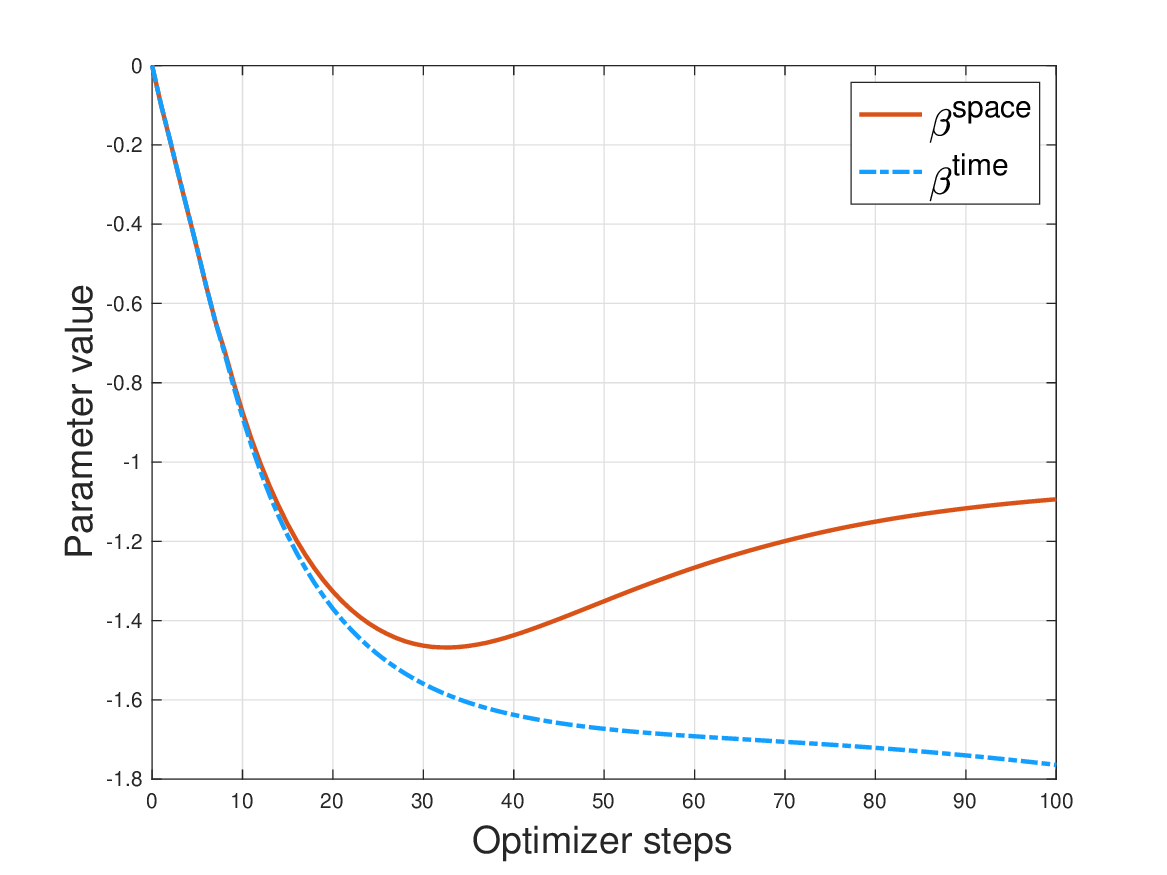}}
	\subfloat[]{\includegraphics[width=0.49\textwidth]{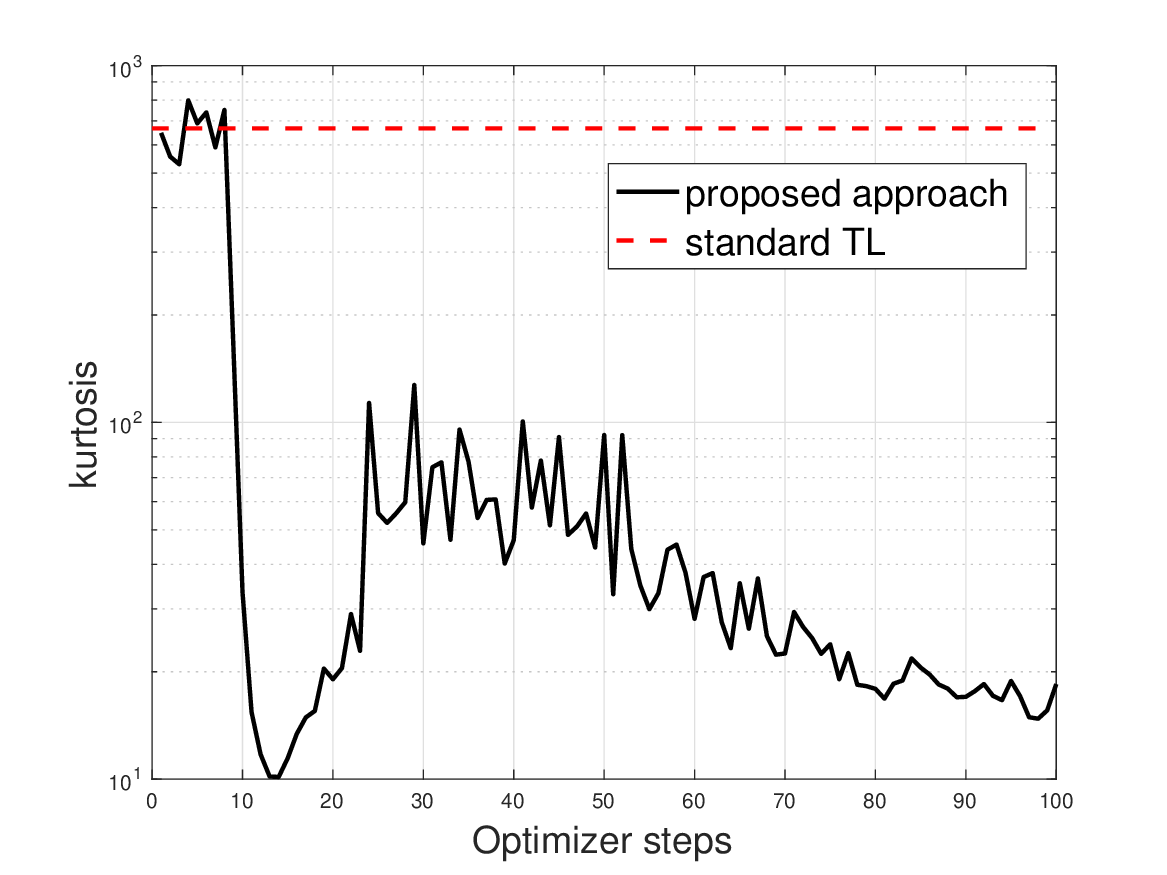}}
\end{figure}

Figure~\ref{fig:decayadam} shows 100 Adam optimization steps for the decay example (Example \ref{exp:decay}) for step size $\Delta t_{pl}=\Delta t_f=1/2^4$. The quantity of interest  is a rare event probability with magnitude $10^{-3}$. To estimate the gradient, we use $M=10^4$ samples per Adam iteration. The squared coefficient of variation is reduced by a factor of $10^{2}$ compared with the standard MC-TL variance after 13 Adam iterations. After reaching this minimum, the squared coefficient of variation increases for the next iteration steps. This behavior might be avoided by employing a smaller step size in the Adam algorithm. Figure~\ref{fig:decayadam}(d) confirms that the kurtosis is bounded to a level below the standard TL's kurtosis, indicating a robust variance estimator.

\begin{figure}
	\caption{Example \ref{exp:mm} with step size $\Delta t_{ pl}=\Delta t_f=1/2^4$ for the proposed IS-MC estimator: \textbf{(a)} sample mean;  \textbf{(b)} squared coefficient of variation; \textbf{(c)} parameters; \textbf{(d)} kurtosis for each optimizer step. The gradient for the Adam optimization, the sample variance, and the kurtosis  were estimated using $M_0=10^5$ samples. Standard MC-TL with step size $\Delta t=1/2^4$ and $M=10^7$ samples was used for comparison.}
	\label{fig:4d}
	\subfloat[]{\includegraphics[width=0.49\textwidth]{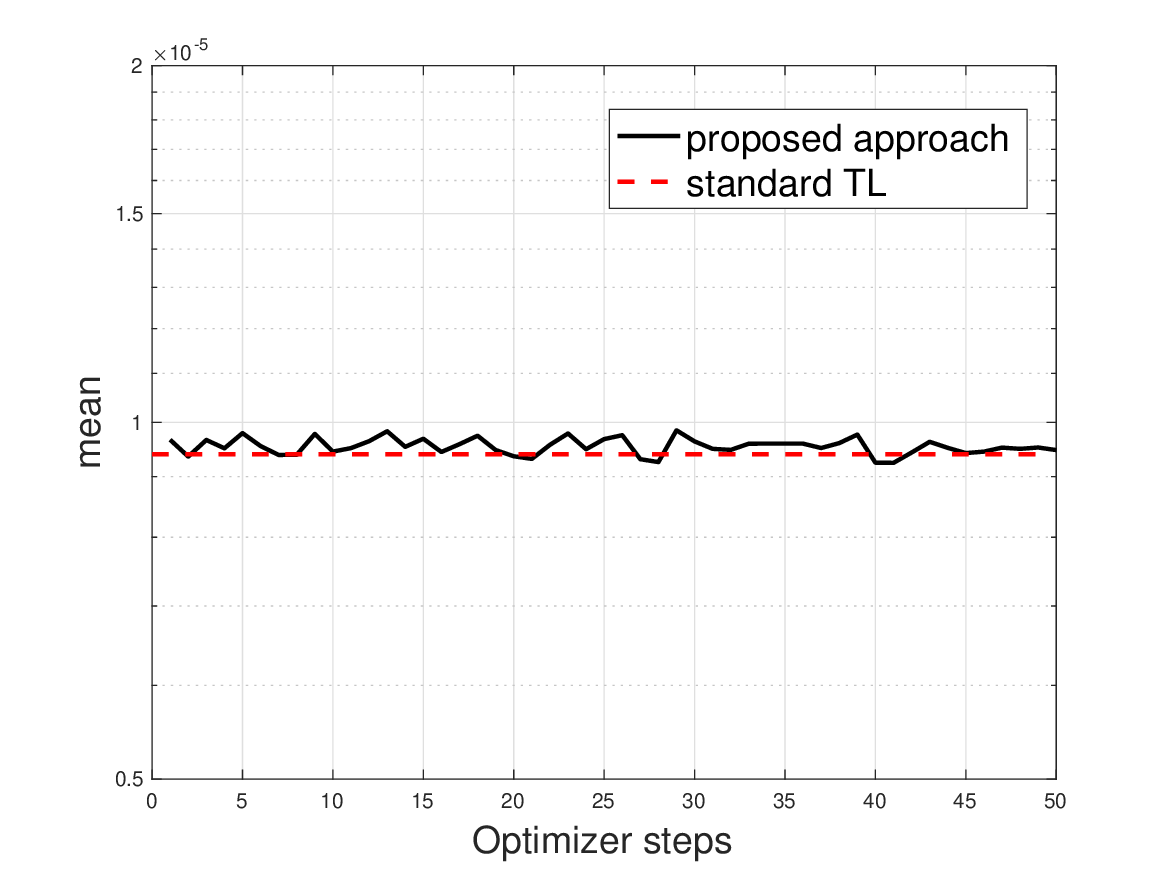}}
	\subfloat[]{\includegraphics[width=0.49\textwidth]{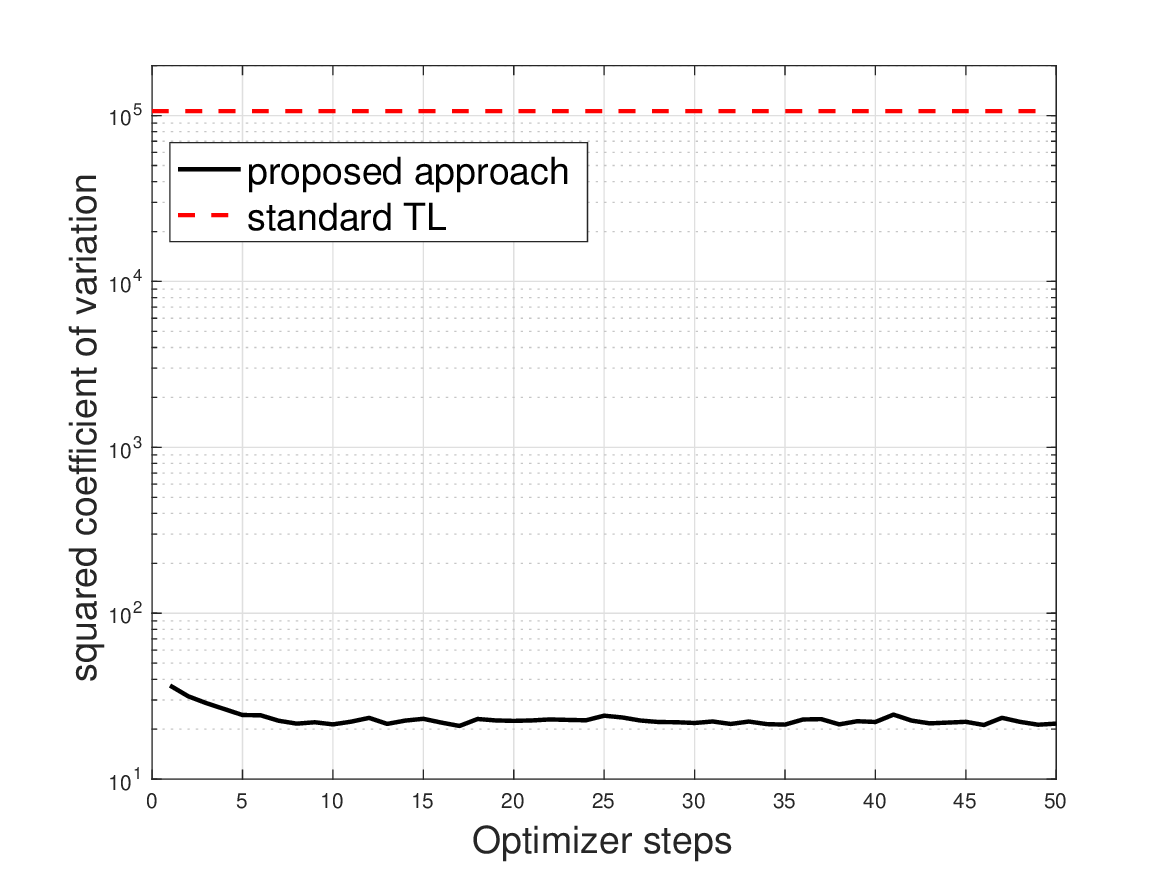}}\\
	\subfloat[ ]{\includegraphics[width=0.49\textwidth]{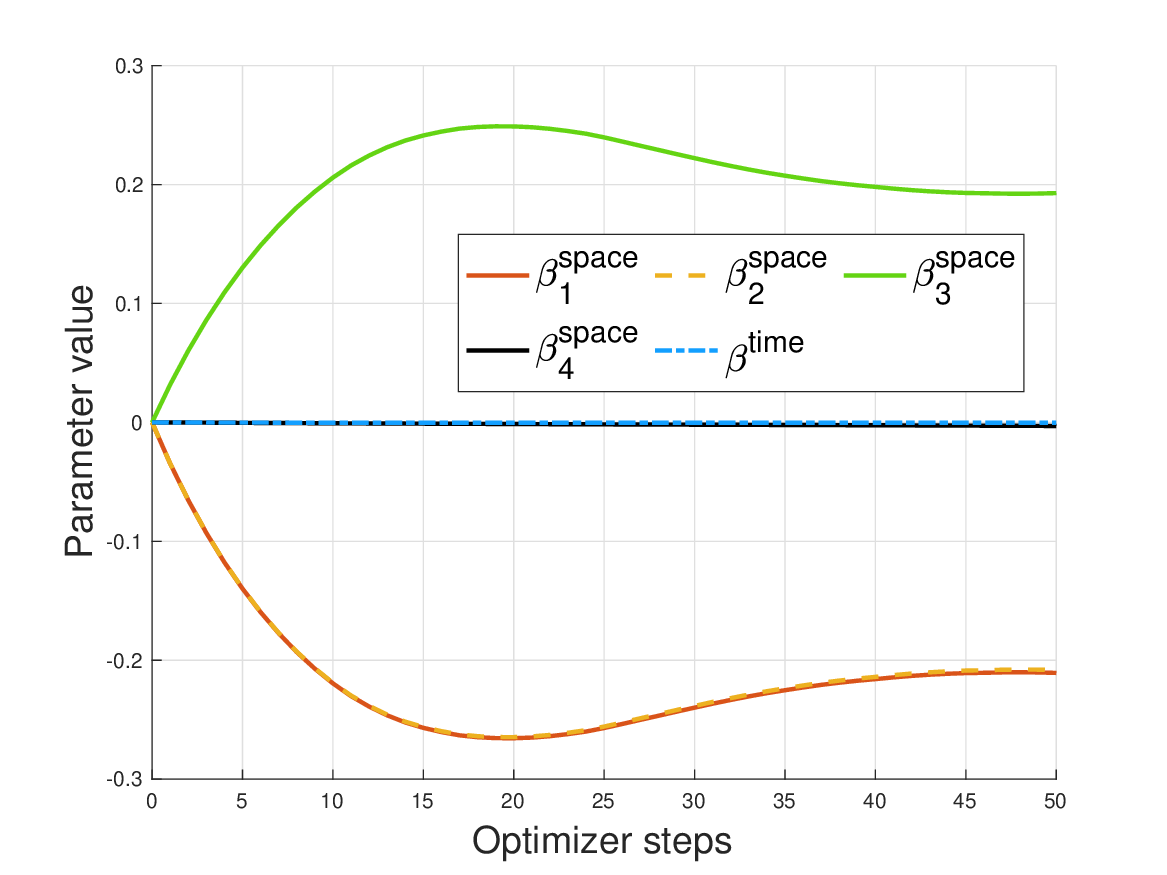}}
	\subfloat[]{\includegraphics[width=0.49\textwidth]{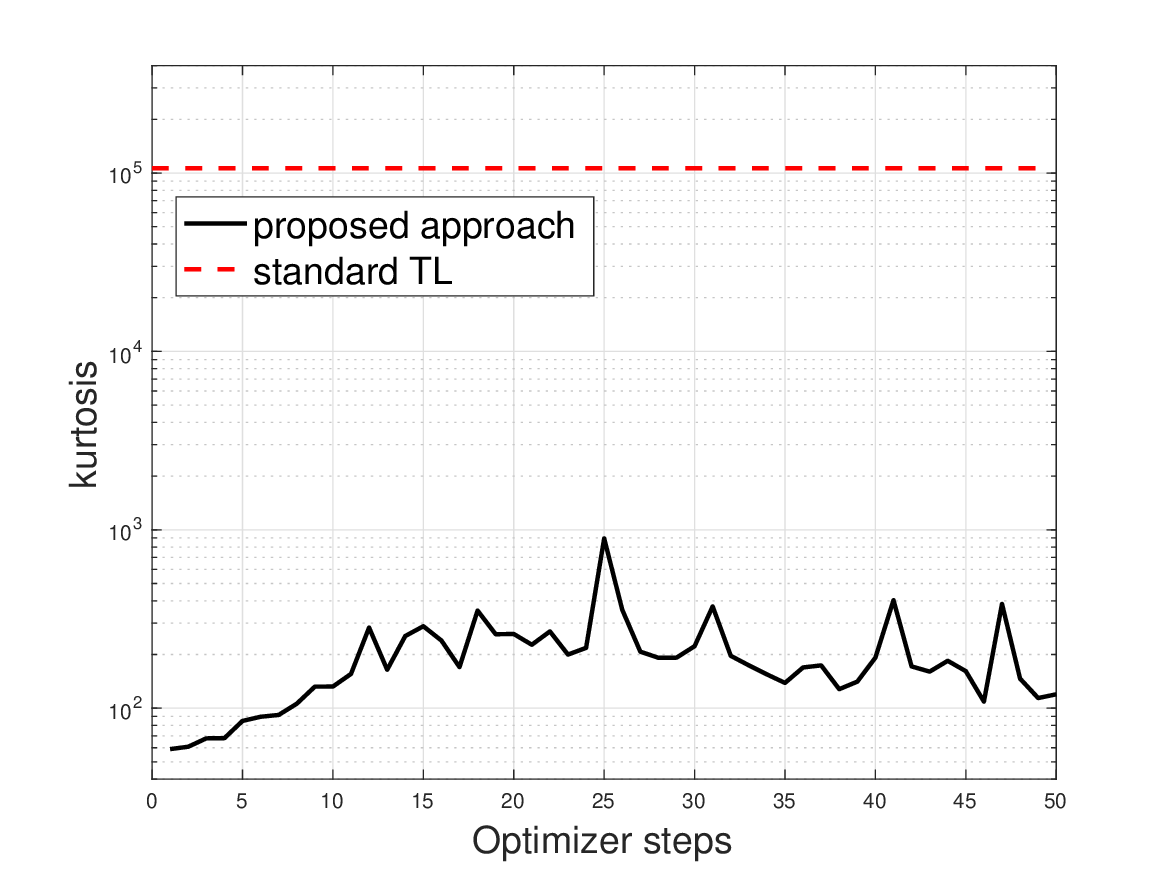}}
\end{figure}

For the 4-dimensional stochastic reaction network (Example~\ref{exp:mm}), the rare event probability for the event $\{X_3(T)>22\}$ is of magnitude $10^{-5}$. Figure~\ref{fig:4d}(b) confirms that the proposed learning-based approach reduces the variance  by a factor $ 4\times 10^3$ compared with standard TL for step size $\Delta t_{pl}=\Delta t_f=1/2^4$. Although Figure~\ref{fig:4d}(c) seems to shows that parameters $\beta^{time}$ and $\beta^{space}_4$ overlap, this is an artifact from the scale of the y-axis; in fact, the final values are $\beta_{time}=-3.2\times 10^{-4}$ and $\beta^{space}_4=-3.0\times 10^{-3}$. The intrinsic structure of Example \ref{exp:mm} results in similar molecule counts for $E(t)$ and $S(t)$ and hence similar values for $\beta^{space}_1$ and $\beta^{space}_2$. Figure~\ref{fig:4d}(d) confirms that the kurtosis for the proposed approach is substantially reduced compared with the kurtosis for the standard TL approach.

\begin{figure}
	\caption{Example~\ref{exp:6d} with step size $\Delta t_{ pl}=\Delta t_f=1/2^4$ for the proposed IS-MC estimators:  \textbf{(a)} sample mean;  \textbf{(b)} squared coefficient of variation; \textbf{(c)} parameters; \textbf{(d)} kurtosis for each optimizer step. The gradient for Adam optimization, the sample variance, and the kurtosis  were estimated using $M_0=10^5$ samples. Standard MC-TL with $M=10^6$ samples and step size $\Delta t=1/2^4$ was used for comparison. 
		\label{fig:6d}}
	\subfloat[]{\includegraphics[width=0.49\textwidth]{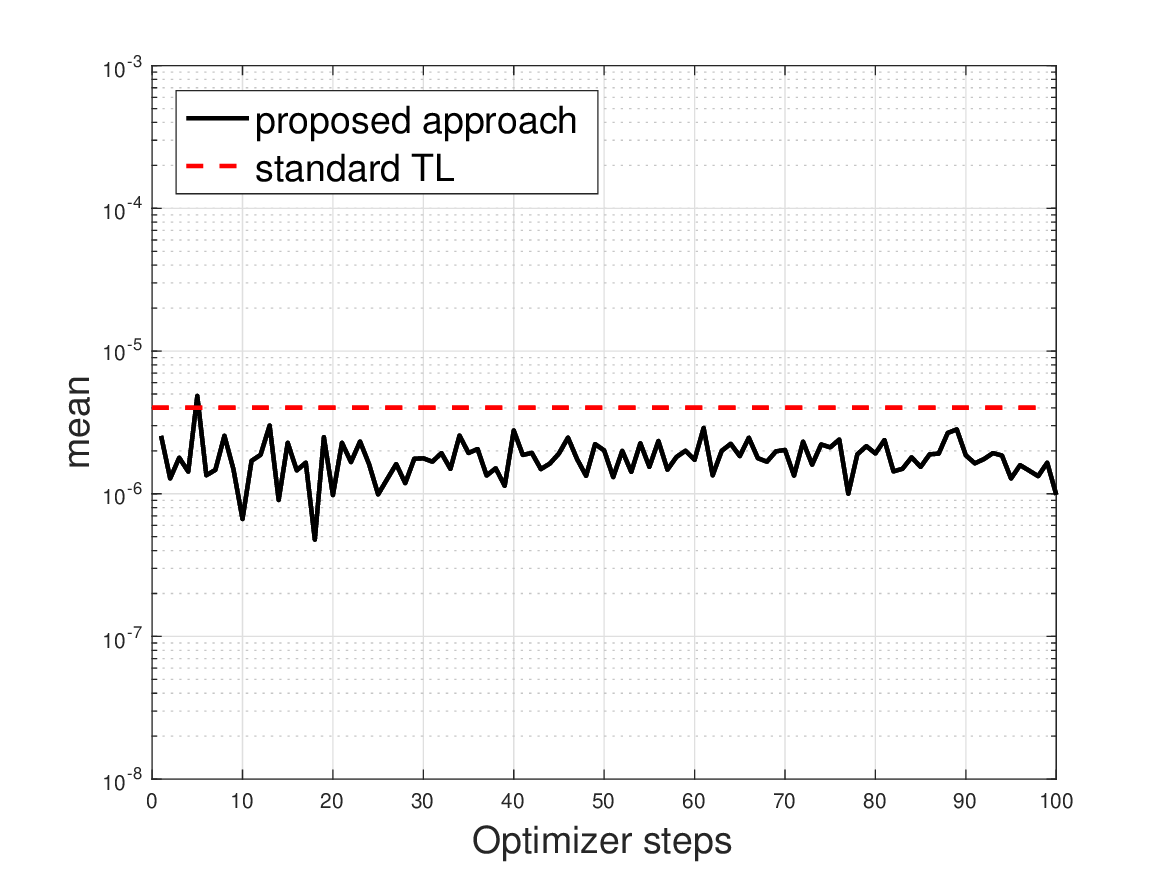}}
	\subfloat[]{\includegraphics[width=0.49\textwidth]{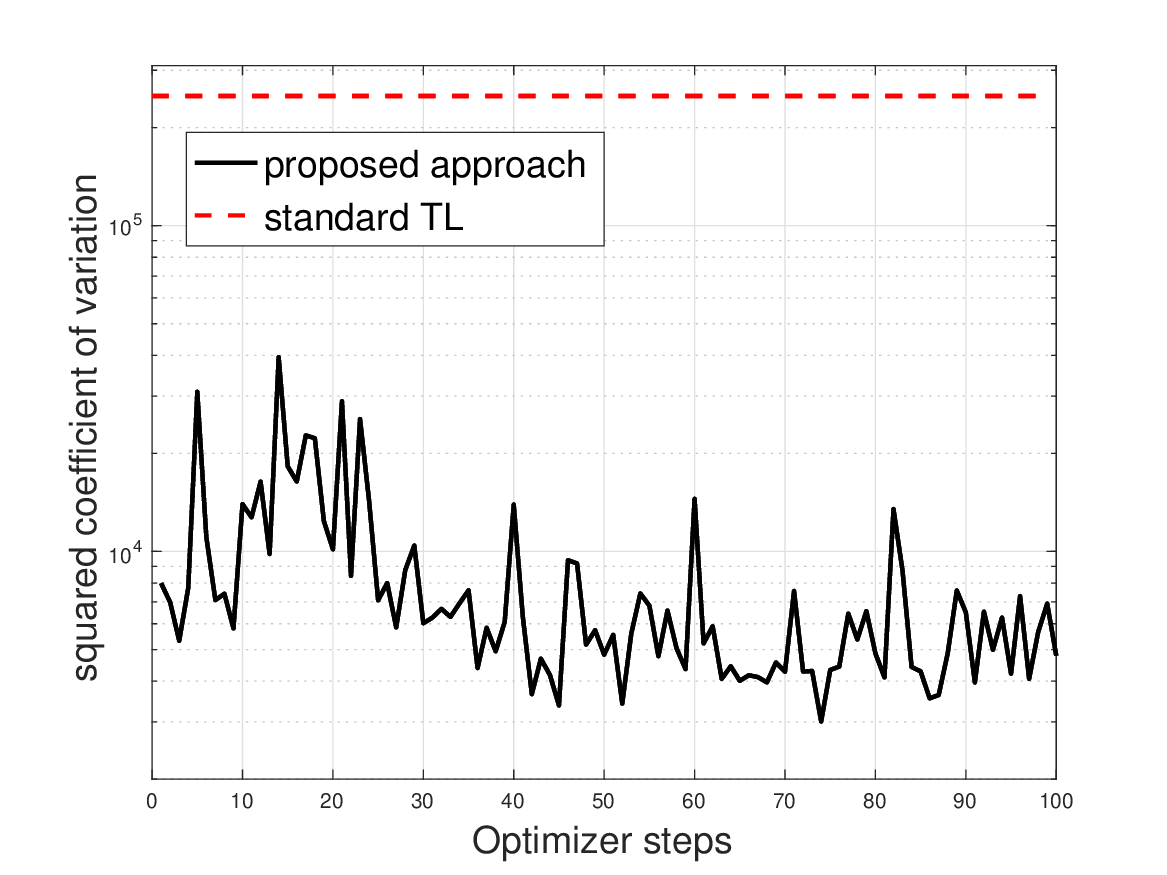}}\\
	\subfloat[ ]{\includegraphics[width=0.49\textwidth]{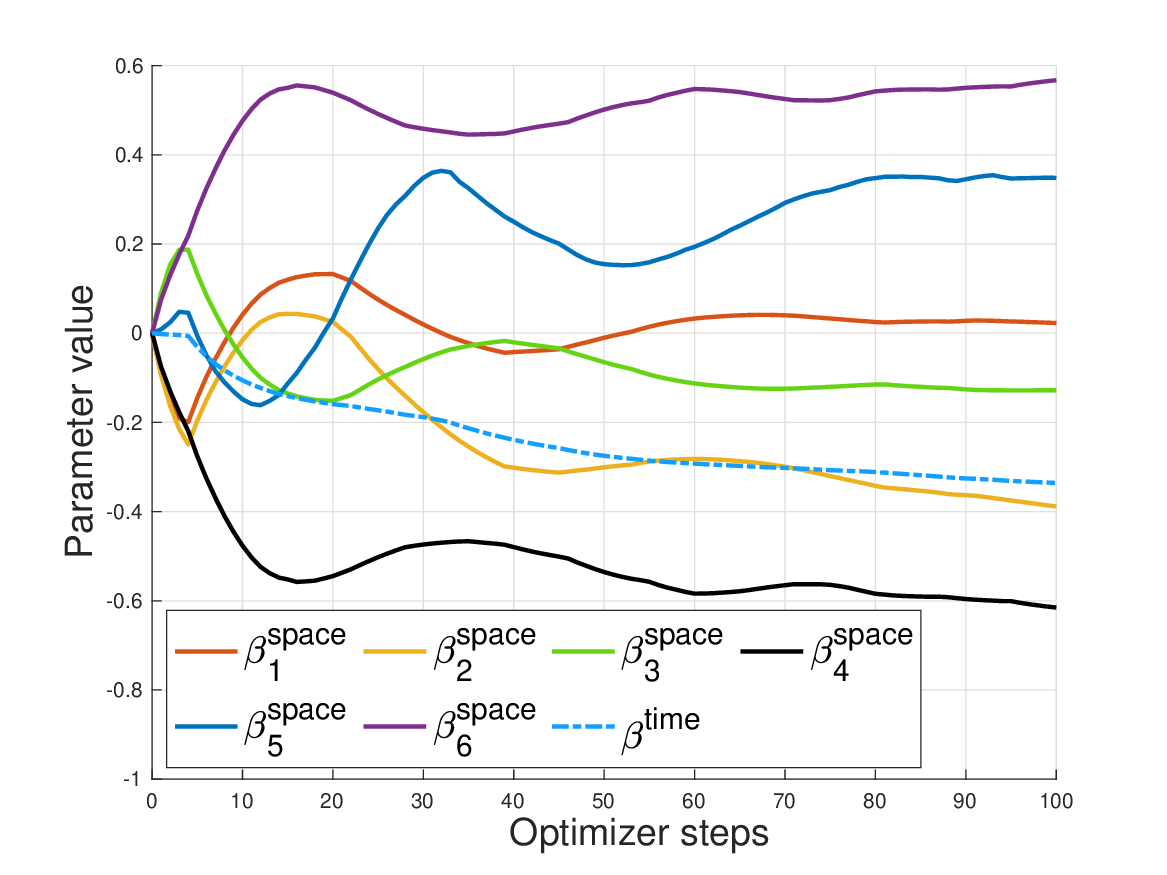}}
	\subfloat[]{\includegraphics[width=0.49\textwidth]{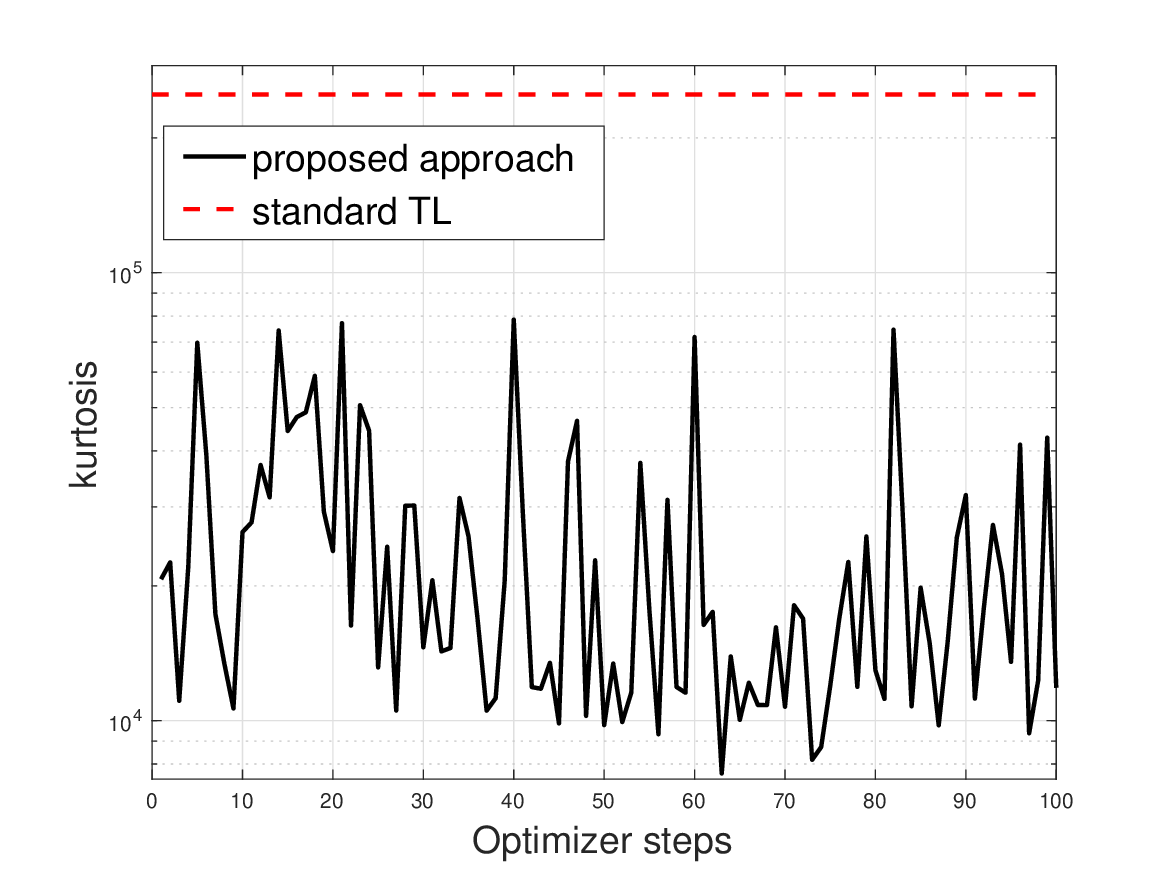}}
\end{figure}

The 6-dimensional example (Example~\ref{exp:6d}) has a rare event probability with magnitude $10^{-6}$. Figure~\ref{fig:6d} shows the Adam optimization results for step size $\Delta t_{pl}=\Delta t_f=1/2^4$. The TL mean differs from the mean for the proposed approach (Fig.~\ref{fig:6d}(a)) because the standard MC-TL estimator requires more than $10^6$ runs to accurately estimate a probability of order $10^{-6}$. The proposed learning-based approach reduces the variance by a factor of  more than $50$ after 43 iterations. The kurtosis is bounded and lower than the kurtosis for the TL approach, confirming that the proposed approach results in a robust variance estimator.

Examples \ref{exp:mm} and \ref{exp:6d} show that a good choice of the ansatz function in combination with reasonable initial parameters provides substantial variance reduction from the first optimization step. However, we do not expect this behavior in general, particularly for high dimensions, and therefore we performed some optimization iterations.

\begin{figure}
	\caption{Example~\ref{exp:mm}, parameters $\boldsymbol{\beta}^{space}$ and $\beta^{time}$ learned with $\Delta t_{pl} =1/2^4$ (see final optimizer step in Figure~\ref{fig:4d}) and applied to forward runs with different $\Delta t_f$ values. The squared coefficient of variation was estimated with $M=10^6$ sample paths. The standard MC-TL approach is used as reference (\textit{dashed red line}).}
	\label{fig:4dII}
	\begin{center}
		\includegraphics[width=0.49\textwidth]{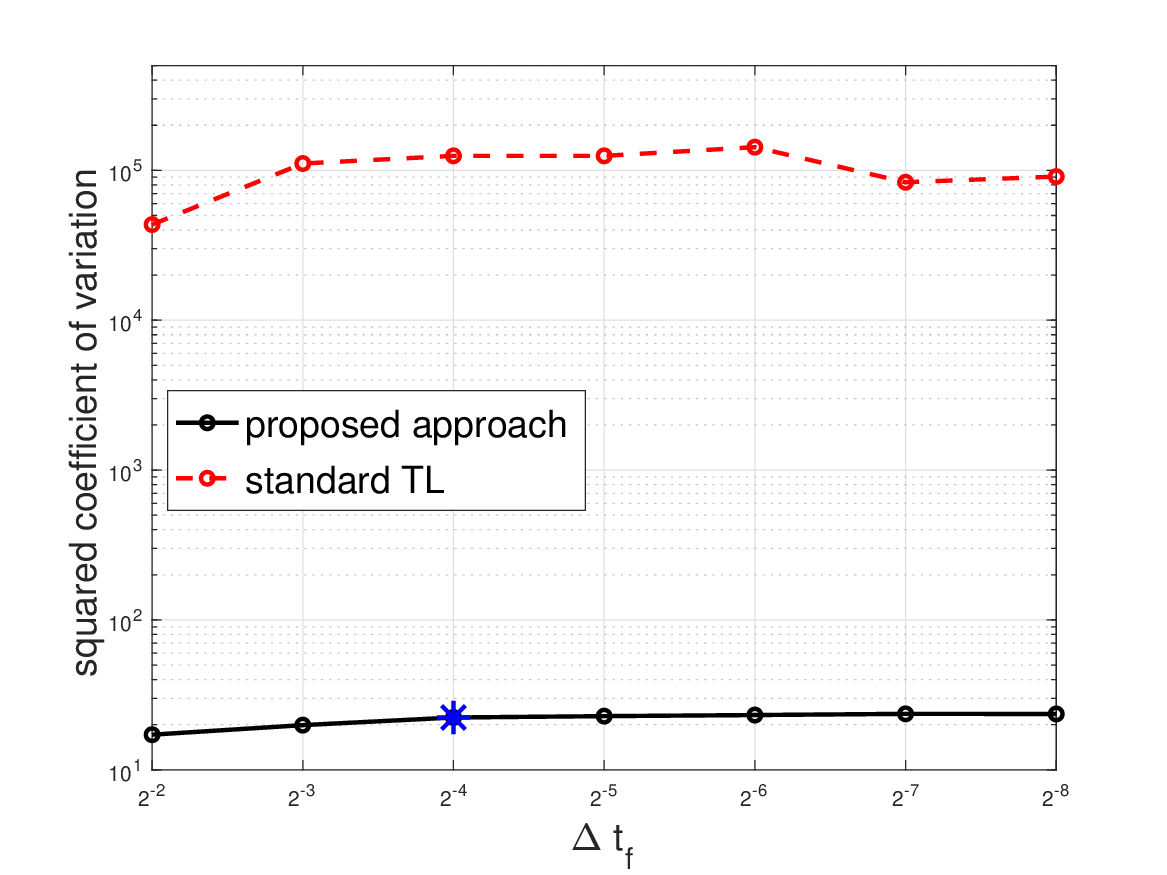}
	\end{center}
\end{figure}

The examples used step size $\Delta t_{pl}=1/2^4$ and showed the squared coefficient of variation with respect to the same step size. To demonstrate that the learned parameters $\boldsymbol{\beta}$ can be used for forward runs with smaller step sizes (\ie,  $\Delta t_{f}\ll\Delta t_{pl}$) as claimed in Section \ref{sec:cost}, we consider Example \ref{exp:mm} and the final parameters from Figure \ref{fig:4dII} for forward runs with different $\Delta t_{f}$. The results show that the variance reduction is constant with respect to $\Delta t_f$, suggesting that a coarse $\Delta t_{pl}$ is sufficient for parameter learning. The same behavior was observed for other tested examples.

\begin{remark}
	We used the ansatz  \eqref{eq:ansatz_sigmoid} based on a single sigmoid for the numerical experiments to demonstrate the potential for the proposed learning-based IS.  Further variance reduction may be achieved either by summing several sigmoid functions as ansatz or selecting  a different basis function shape. Relevant analyses will be pursued in future work. 
\end{remark}

\section{Conclusions and Future Work}
\label{sec:Conclusions and future work}

This work developed an efficient path-dependent IS scheme to estimate statistical quantities for SRN processes, particularly rare event probabilities.  Optimal IS parameters were obtained within a pre-selected class of change of measure using the proposed connection to an associated SOC problem, which could be solved via dynamic programming. To mitigate the curse of dimensionality encountered by the dynamic programming relation, we proposed a method for multi-dimensional SRNs based on approximating the value function via an ansatz function (\ie~a neural network), where the parameters were learned using a stochastic optimization algorithm. Numerical examples and subsequent analyses verified that the proposed estimator achieved substantial variance reduction compared with the standard MC method, providing lowered computational complexity in the rare event regime.

Future work will further analyze the proposed learning-based approach and expand it to derive a multilevel MC estimator. We also plan to combine an implementation of the dynamic programming principle as derived in Section \ref{sec:Algorithm} with dimension reduction methods for SRNs.

\textbf{Acknowledgments}
This publication is based upon work supported by the King Abdullah University of Science and Technology (KAUST) Office of Sponsored Research (OSR) under Award No. OSR-2019-CRG8-4033.
This work was partially performed as part of the Helmholtz School for Data Science in Life, Earth and Energy (HDS-LEE) and received funding from the Helmholtz Association of German Research Centres and the Alexander von Humboldt Foundation. 

\bibliographystyle{plain}
\bibliography{IS_via_SOC_manuscript} 

\begin{thebibliography}{10}

\bibitem{abdulle2010chebyshev}
Assyr Abdulle, Yucheng Hu, and Tiejun Li.
\newblock Chebyshev methods with discrete noise: the $\tau$-{ROCK} methods.
\newblock {\em Journal of Computational Mathematics}, pages 195--217, 2010.

\bibitem{ahn2013implicit}
Tae-Hyuk Ahn, Adrian Sandu, and Xiaoying Han.
\newblock Implicit simulation methods for stochastic chemical kinetics.
\newblock {\em arXiv preprint arXiv:1303.3614}, 2013.

\bibitem{Anderson2012}
D.~{Anderson} and D.~{Higham}.
\newblock Multilevel {M}onte {C}arlo for continuous {M}arkov chains, with
  applications in biochemical kinetics.
\newblock {\em SIAM Multiscal Model. Simul.}, 10(1), 2012.

\bibitem{anderson2007modified}
David~F Anderson.
\newblock A modified next reaction method for simulating chemical systems with
  time dependent propensities and delays.
\newblock {\em The Journal of chemical physics}, 127(21):214107, 2007.

\bibitem{anderson2015stochastic}
David~F Anderson and Thomas~G Kurtz.
\newblock {\em Stochastic analysis of biochemical systems}, volume~1.
\newblock Springer, 2015.

\bibitem{aparicio2001population}
Juan~P Aparicio and Hern{\'a}n~G Solari.
\newblock Population dynamics: Poisson approximation and its relation to the
  langevin process.
\newblock {\em Physical Review Letters}, 86(18):4183, 2001.

\bibitem{banisch2013meshless}
Ralf Banisch and Carsten Hartmann.
\newblock Meshless discretization of lq-type stochastic control problems.
\newblock {\em arXiv preprint arXiv:1309.7497}, 2013.

\bibitem{bayer2016efficient}
Christian Bayer, Alvaro Moraes, Ra{\'u}l Tempone, and Pedro Vilanova.
\newblock An efficient forward--reverse expectation-maximization algorithm for
  statistical inference in stochastic reaction networks.
\newblock {\em Stochastic Analysis and Applications}, 34(2):193--231, 2016.

\bibitem{ben2020hierarchical}
Chiheb Ben~Hammouda.
\newblock {\em Hierarchical Approximation Methods for Option Pricing and
  Stochastic Reaction Networks}.
\newblock PhD thesis, 2020.

\bibitem{hammouda2020importance}
Chiheb Ben~Hammouda, Nadhir Ben~Rached, and Ra{\'u}l Tempone.
\newblock Importance sampling for a robust and efficient multilevel {M}onte
  {C}arlo estimator for stochastic reaction networks.
\newblock {\em Statistics and Computing}, 30(6):1665--1689, 2020.

\bibitem{hammouda2017multilevel}
Chiheb Ben~Hammouda, Alvaro Moraes, and Ra{\'u}l Tempone.
\newblock Multilevel hybrid split-step implicit tau-leap.
\newblock {\em Numerical Algorithms}, 74(2):527--560, 2017.

\bibitem{ben2021efficient}
Nadhir Ben~Rached, Abdul-Lateef Haji-Ali, Gerardo Rubino, and Ra{\'u}l Tempone.
\newblock Efficient importance sampling for large sums of independent and
  identically distributed random variables.
\newblock {\em Statistics and Computing}, 31(6):1--13, 2021.

\bibitem{brauer2001mathematical}
Fred Brauer and Carlos Castillo-Chavez.
\newblock {\em Mathematical models in population biology and epidemiology},
  volume~40.
\newblock Springer.

\bibitem{cao2005trapezoidal}
Yang Cao and Linda Petzold.
\newblock Trapezoidal tau-leaping formula for the stochastic simulation of
  biochemical systems.
\newblock {\em Proceedings of Foundations of Systems Biology in Engineering
  (FOSBE 2005)}, pages 149--152, 2005.

\bibitem{cao2013adaptively}
Youfang Cao and Jie Liang.
\newblock Adaptively biased sequential importance sampling for rare events in
  reaction networks with comparison to exact solutions from finite buffer
  d{CME} method.
\newblock {\em The Journal of chemical physics}, 139(2):07B605\_1, 2013.

\bibitem{daigle2011automated}
Bernie~J Daigle~Jr, Min~K Roh, Dan~T Gillespie, and Linda~R Petzold.
\newblock Automated estimation of rare event probabilities in biochemical
  systems.
\newblock {\em The Journal of chemical physics}, 134(4):01B628, 2011.

\bibitem{dupuis2012importance}
Paul Dupuis, Konstantinos Spiliopoulos, and Hui Wang.
\newblock Importance sampling for multiscale diffusions.
\newblock {\em Multiscale Modeling \& Simulation}, 10(1):1--27, 2012.

\bibitem{engblom2012stability}
Stefan Engblom.
\newblock On the stability of stochastic jump kinetics.
\newblock {\em arXiv preprint arXiv:1202.3892}, 2012.

\bibitem{kurtz_2005}
Stewart~N. Ethier and Thomas~G. Kurtz.
\newblock {\em Markov processes : characterization and convergence}.
\newblock Wiley series in probability and mathematical statistics. J. Wiley \&
  Sons, New York, Chichester, 1986.

\bibitem{fleming2006controlled}
Wendell~H Fleming and Halil~Mete Soner.
\newblock {\em Controlled Markov processes and viscosity solutions}, volume~25.
\newblock Springer Science \& Business Media, 2006.

\bibitem{giles2008multilevel}
Michael~B Giles.
\newblock Multilevel {M}onte {C}arlo path simulation.
\newblock {\em Operations Research}, 56(3):607--617, 2008.

\bibitem{giles2015multilevel}
Michael~B Giles.
\newblock Multilevel {M}onte {C}arlo methods.
\newblock {\em Acta Numerica}, 24:259--328, 2015.

\bibitem{gillespie2019guided}
Colin~S Gillespie and Andrew Golightly.
\newblock Guided proposals for efficient weighted stochastic simulation.
\newblock {\em The Journal of chemical physics}, 150(22):224103, 2019.

\bibitem{gillespie2009refining}
Dan~T Gillespie, Min Roh, and Linda~R Petzold.
\newblock Refining the weighted stochastic simulation algorithm.
\newblock {\em The Journal of chemical physics}, 130(17):174103, 2009.

\bibitem{gillespie1976general}
Daniel~T Gillespie.
\newblock A general method for numerically simulating the stochastic time
  evolution of coupled chemical reactions.
\newblock {\em Journal of computational physics}, 22(4):403--434, 1976.

\bibitem{gillespie2001approximate}
Daniel~T Gillespie.
\newblock Approximate accelerated stochastic simulation of chemically reacting
  systems.
\newblock {\em The Journal of Chemical Physics}, 115(4):1716--1733, 2001.

\bibitem{gupta2014scalable}
Ankit Gupta, Corentin Briat, and Mustafa Khammash.
\newblock A scalable computational framework for establishing long-term
  behavior of stochastic reaction networks.
\newblock {\em PLoS computational biology}, 10(6):e1003669, 2014.

\bibitem{hartmann2014characterization}
Carsten Hartmann, Ralf Banisch, Marco Sarich, Tomasz Badowski, and Christof
  Sch{\"u}tte.
\newblock Characterization of rare events in molecular dynamics.
\newblock {\em Entropy}, 16(1):350--376, 2014.

\bibitem{hartmann2019variational}
Carsten Hartmann, Omar Kebiri, Lara Neureither, and Lorenz Richter.
\newblock Variational approach to rare event simulation using least-squares
  regression.
\newblock {\em Chaos: An Interdisciplinary Journal of Nonlinear Science},
  29(6):063107, 2019.

\bibitem{hartmann2017variational}
Carsten Hartmann, Lorenz Richter, Christof Sch{\"u}tte, and Wei Zhang.
\newblock Variational characterization of free energy: Theory and algorithms.
\newblock {\em Entropy}, 19(11):626, 2017.

\bibitem{hartmann2018importance}
Carsten Hartmann, Christof Sch{\"u}tte, Marcus Weber, and Wei Zhang.
\newblock Importance sampling in path space for diffusion processes with
  slow-fast variables.
\newblock {\em Probability Theory and Related Fields}, 170(1):177--228, 2018.

\bibitem{hensel2009stochastic}
Sebastian~C Hensel, James~B Rawlings, and John Yin.
\newblock Stochastic kinetic modeling of vesicular stomatitis virus
  intracellular growth.
\newblock {\em Bulletin of mathematical biology}, 71(7):1671--1692, 2009.

\bibitem{kebiri2017adaptive}
Omar Kebiri, Lara Neureither, and Carsten Hartmann.
\newblock Adaptive importance sampling with forward-backward stochastic
  differential equations.
\newblock In {\em International workshop on Stochastic Dynamics out of
  Equilibrium}, pages 265--281. Springer, 2017.

\bibitem{kingma2014adam}
Diederik~P Kingma and Jimmy Ba.
\newblock Adam: A method for stochastic optimization.
\newblock {\em arXiv preprint arXiv:1412.6980}, 2014.

\bibitem{kroese2013handbook}
Dirk~P Kroese, Thomas Taimre, and Zdravko~I Botev.
\newblock {\em Handbook of {M}onte {C}arlo methods}.
\newblock John Wiley \& Sons, 2013.

\bibitem{kuwahara2008efficient}
Hiroyuki Kuwahara and Ivan Mura.
\newblock An efficient and exact stochastic simulation method to analyze rare
  events in biochemical systems.
\newblock {\em The Journal of chemical physics}, 129(16):10B619, 2008.

\bibitem{l1995note}
Pierre L'Ecuyer.
\newblock Note: On the interchange of derivative and expectation for likelihood
  ratio derivative estimators.
\newblock {\em Management Science}, 41(4):738--747, 1995.

\bibitem{lester2015adaptive}
Christopher Lester, Christian~Adam Yates, Michael~B Giles, and Ruth~E Baker.
\newblock An adaptive multi-level simulation algorithm for stochastic
  biological systems.
\newblock {\em The Journal of chemical physics}, 142(2):01B612\_1, 2015.

\bibitem{li2007analysis}
Tiejun Li.
\newblock Analysis of explicit tau-leaping schemes for simulating chemically
  reacting systems.
\newblock {\em Multiscale Modeling \& Simulation}, 6(2):417--436, 2007.

\bibitem{moraes2016multilevel_splitting}
Alvaro Moraes, Ra{\'u}l Tempone, and Pedro Vilanova.
\newblock A multilevel adaptive reaction-splitting simulation method for
  stochastic reaction networks.
\newblock {\em SIAM Journal on Scientific Computing}, 38(4):A2091--A2117, 2016.

\bibitem{moraes2016multilevel}
Alvaro Moraes, Ra{\'u}l Tempone, and Pedro Vilanova.
\newblock A multilevel adaptive reaction-splitting simulation method for
  stochastic reaction networks.
\newblock {\em SIAM Journal on Scientific Computing}, 38(4):A2091--A2117, 2016.

\bibitem{nusken2021solving}
Nikolas N{\"u}sken and Lorenz Richter.
\newblock Solving high-dimensional hamilton--jacobi--bellman pdes using neural
  networks: perspectives from the theory of controlled diffusions and measures
  on path space.
\newblock {\em Partial Differential Equations and Applications}, 2(4):1--48,
  2021.

\bibitem{rached2022double}
Nadhir~Ben Rached, Abdul-Lateef Haji-Ali, Shyam Mohan, and Ra{\'u}l Tempone.
\newblock Double loop {M}onte {C}arlo estimator with importance sampling for
  mckean-vlasov stochastic differential equation.
\newblock {\em arXiv preprint arXiv:2207.06926}, 2022.

\bibitem{rao2003stochastic}
Christopher~V Rao and Adam~P Arkin.
\newblock Stochastic chemical kinetics and the quasi-steady-state assumption:
  Application to the {G}illespie algorithm.
\newblock {\em The Journal of chemical physics}, 118(11):4999--5010, 2003.

\bibitem{rathinam2013moment}
Muruhan Rathinam.
\newblock Moment growth bounds on continuous time {M}arkov processes on
  non-negative integer lattices.
\newblock {\em arXiv preprint arXiv:1304.5169}, 2013.

\bibitem{rathinam2007reversible}
Muruhan Rathinam and Hana El~Samad.
\newblock Reversible-equivalent-monomolecular tau: A leaping method for
  “small number and stiff” stochastic chemical systems.
\newblock {\em Journal of Computational Physics}, 224(2):897--923, 2007.

\bibitem{roh2019data}
Min~K Roh.
\newblock Data-driven method for efficient characterization of rare event
  probabilities in biochemical systems.
\newblock {\em Bulletin of mathematical biology}, 81(8):3097--3120, 2019.

\bibitem{roh2010state}
Min~K Roh, Dan~T Gillespie, and Linda~R Petzold.
\newblock State-dependent biasing method for importance sampling in the
  weighted stochastic simulation algorithm.
\newblock {\em The Journal of chemical physics}, 133(17):174106, 2010.

\bibitem{srivastava2002stochastic}
Ranjan Srivastava, L~You, J~Summers, and J~Yin.
\newblock Stochastic vs. deterministic modeling of intracellular viral
  kinetics.
\newblock {\em Journal of theoretical biology}, 218(3):309--321, 2002.

\bibitem{zhang2014applications}
Wei Zhang, Han Wang, Carsten Hartmann, Marcus Weber, and Christof Sch{\"u}tte.
\newblock Applications of the cross-entropy method to importance sampling and
  optimal control of diffusions.
\newblock {\em SIAM Journal on Scientific Computing}, 36(6):A2654--A2672, 2014.

\end{thebibliography}
\pagebreak
\appendix
\section{Proof for Theorem~\ref{theo:exact_optival}}
\label{appendix:Proof of Theorem dynamic prog}
\begin{proof}[Proof for Theorem~\ref{theo:exact_optival}]
	To show (\ref{eq:exact_optival}), we first reformulate $C_{n,x}(\boldsymbol{\delta}^{\Delta t}_n,\dots,\boldsymbol{\delta}^{\Delta t}_{N-1})$ using the definition for the likelihood and the notion of conditional expectation,
	\begin{small}
		\begin{align}\label{eq:proof1}
			C_{n,x}&(\boldsymbol{\delta}^{\Delta t}_n,\dots,\boldsymbol{\delta}^{\Delta t}_{N-1})\nonumber\\
			&=\mathbb{E}\left[g^2(\overline{X}_N^{\Delta t})\prod_{k=n}^{N-1} L_k^2(\bar{\mathbf{P}}_k,\boldsymbol{\delta}_k^{\Delta t}(\overline{X}_k^{\Delta t})) \mid \overline{X}_n^{\Delta t}=x\right]\nonumber
			\\
			&=\mathbb{E}\left[g^2(\overline{X}_N^{\Delta t})\cdot  L_n^2(\bar{\mathbf{P}}_n,\boldsymbol{\delta}_n^{\Delta t}(\overline{X}^{\Delta t}_n)) \cdot \prod_{k=n+1}^{N-1} L_k^2(\bar{\mathbf{P}}_k,\boldsymbol{\delta}_k^{\Delta t}(\overline{X}_k^{\Delta t})) \mid \overline{X}_n^{\Delta t}=x\right]\nonumber\\
			&=\mathbb{E}\left[
			g^2(\overline{X}_N^{\Delta t})\cdot \exp\left(- 2 \left(\sum_{j=1}^J a_j(\overline{X}_{n}^{\Delta t})-\delta_{n,j}^{\Delta t}(\overline{X}^{\Delta t}_n)\right)\Delta t\right)\left(\prod_{j=1}^J\frac{a_j(\overline{X}_{n}^{\Delta t})}{\delta_{n,j}^{\Delta t}(\overline{X}^{\Delta t}_n)}\right)^{2\bar{P}_{n,j}} \right.\nonumber\\
			& \left.\cdot \prod_{k=n+1}^{N-1} L_k^2(\bar{\mathbf{P}}_k,\boldsymbol{\delta}_k^{\Delta t}(\overline{X}_k^{\Delta t})) \mid \overline{X}_n^{\Delta t}=x\right].
		\end{align}
		
		Setting
		\begin{small}
			\begin{equation*}
				B(\bar{\mathbf{P}}_n):=g^2(\overline{X}_N^{\Delta t}) \exp\left(-2\left(\sum_{j=1}^J a_j(\overline{X}_{n}^{\Delta t})-\delta_{n,j}^{\Delta t}(\overline{X}^{\Delta t}_n)\right)\Delta t\right)\left(\prod_{j=1}^J\frac{a_j(\overline{X}_{n}^{\Delta t})}{\delta_{n,j}^{\Delta t}(\overline{X}^{\Delta t}_n)}\right)^{2\bar{P}_{n,j}} \prod_{k=n+1}^{N-1} L_k^2(\bar{\mathbf{P}}_k,\boldsymbol{\delta}_k^{\Delta t}(\overline{X}_k^{\Delta t})),
			\end{equation*}
		\end{small}
		we can reformulate \eqref{eq:proof1} and derive
		\begin{align*}
			C_{n,x}&(\boldsymbol{\delta}^{\Delta t}_n,\dots,\boldsymbol{\delta}^{\Delta t}_{N-1})\nonumber\\
			&= \sum_{\mathbf{p}\in \mathbb{N}^J}\mathbb{P}\left(\bar{\mathbf{P}}_n=\mathbf{p} \mid \overline{X}_n^{\Delta t}=x\right) \cdot \mathbb{E} \left[B(\bar{\mathbf{P}}_n) \mid \overline{X}_n^{\Delta t}=x,\bar{\mathbf{P}}_n=\mathbf{p}\right] \\
			&=\sum_{\mathbf{p}\in \mathbb{N}^J} \left[\prod_{j=1}^{J} \frac{(\Delta t \cdot \delta_{n,j}^{\Delta t}(x))^{p_j}}{p_j!}\exp(-\Delta t \cdot \sum_{j=0}^J \delta_{n,j}^{\Delta t}(x))\right] \cdot \exp\left(- 2 \left (\sum_{j=1}^J a_j(x)-\delta_{n,j}^{\Delta t}(x)\right)\Delta t\right)\\
			&~~~~~~~~~~\cdot\left(\prod_{j=1}^J\frac{a_j(x)}{\boldsymbol{\delta}_{n,j}^{\Delta t}(x)}\right)^{2p_j} \cdot \mathbb{E}\left[g^2(\overline{X}_N^{\Delta t})\prod_{k=n+1}^{N-1} L_k^2\left(\bar{\mathbf{P}}_k,\delta_k^{\Delta t}(\overline{X}_k^{\Delta t})\right) \mid \overline{X}_n^{\Delta t}=x,\bar{\mathbf{P}}_n=\mathbf{p}\right]\\
			&=\exp\left(\left(-2\sum_{j=1}^J a_j(x)+\sum_{j=1}^J\delta_{n,j}^{\Delta t}(x)\right)\Delta t\right)\cdot\sum_{\mathbf{p}\in \mathbb{N}^J}\left(\prod_{j=1}^{J} \frac{(\Delta t \cdot \delta_{n,j}^{\Delta t}(x))^{p_j}}{p_j!} (\frac{a_j(x)}{\delta_{n,j}^{\Delta t}(x)})^{2p_j}\right) \\
			&~~~~~~~~~~\cdot \mathbb{E}\left[g^2(\overline{X}_N^{\Delta t})\prod_{k=n+1}^{N-1} L_k^2 \left(\bar{\mathbf{P}}_k,\boldsymbol{\delta}_k^{\Delta t}(\overline{X}_k^{\Delta t})\right) \mid \overline{X}_{n+1}^{\Delta t}=\max(0,x +\mathbf{p}^T \boldsymbol{\nu})\right].
		\end{align*}
	\end{small}
	We can prove Theorem~\ref{theo:exact_optival} using the above results. We split the proof into two  parts, where  the first inequality is obtained by  
	\begin{small}
		\begin{align*}
			&u_{\Delta t}(n,x)\\
			&= \inf_{\{\boldsymbol{\delta}^{\Delta t}_i\}_{i=n,\dots,N-1} \in \mathcal{A}^{N-n}} \Bigg[\exp\left(\left(-2\sum_{j=1}^J a_j(x)+\sum_{j=1}^J\delta_{n,j}^{\Delta t}(x)\right)\Delta t\right) \sum_{\mathbf{p}\in \mathbb{N}^J} \left(\left(\prod_{j=1}^{J} \frac{(\Delta t \cdot \delta_{n,j}^{\Delta t}(x))^{p_j}}{p_j!} (\frac{a_j(x)}{\delta_{n,j}^{\Delta t}(x)})^{2p_j} \right) \right.\\
			&  \left. ~~~~~~~~~~~~~~~~~~~~~~~~~~~ ~~~~~~~~~~~~~~~~~~~~\times \mathbb{E}\left[g^2(\overline{X}_N^{\Delta t})\prod_{k=n+1}^{N-1} L_k^2(\bar{\mathbf{P}}_k,\boldsymbol{\delta}_k^{\Delta t}(\overline{X}_k^{\Delta t})) \mid \overline{X}_{n+1}^{\Delta t}=\max(0,x +\mathbf{p}^T\nu)\right] \right)\Bigg]\\
			&\geq  \inf_{\{\boldsymbol{\delta}^{\Delta t}_i\}_{i=n,\dots,N-1} \in \mathcal{A}^{N-n}}\Bigg[ \exp\left(\left(-2\sum_{j=1}^J a_j(x)+\sum_{j=1}^J\delta_{n,j}^{\Delta t}(x)\right)\Delta t\right) \sum_{\mathbf{p}\in \mathbb{N}^J} \left(\left(\prod_{j=1}^{J} \frac{(\Delta t \cdot \delta_{n,j}^{\Delta t}(x))^{p_j}}{p_j!} (\frac{a_j(x)}{\delta_{n,j}^{\Delta t}(x)})^{2p_j} \right)\right.\\
			&\left. ~~~~~~~~~~~~~~\times \inf_{\{\boldsymbol{\delta}^{\Delta t}_k\}_{k=n+1,\dots,N-1} \in \mathcal{A}^{N-n-1}} \mathbb{E}\left[g^2(\overline{X}_N^{\Delta t})\prod_{k=n+1}^{N-1} L_k^2(\bar{\mathbf{P}}_k,\boldsymbol{\delta}_k^{\Delta t}(\overline{X}_k^{\Delta t}))\mid \overline{X}_{n+1}^{\Delta t}=\max(0,x +\mathbf{p}^T \boldsymbol{\nu})\right] \right) \Bigg]\\
			&= \inf_{\boldsymbol{\delta}^{\Delta t}_n(x) \in \mathcal{A}_x} \Bigg[ \exp\left(\left(-2\sum_{j=1}^J a_j(x)+\sum_{j=1}^J\delta_{n,j}^{\Delta t}(x)\right)\Delta t\right)\\
			&  ~~~~~~~~~~~~~~~~~~~~\times\sum_{\mathbf{p}\in \mathbb{N}^J} \left(\left(\prod_{j=1}^{J} \frac{(\Delta t \cdot \delta_{n,j}^{\Delta t}(x))^{p_j}}{p_j!} \left(\frac{a_j(x)}{\delta_{n,j}^{\Delta t}(x)}\right)^{2p_j} \right)\cdot u_{\Delta t}(n+1,\max(0,x+\mathbf{p}^T \boldsymbol{\nu})) \right)\Bigg]
		\end{align*}
	\end{small}
	To prove the second inequality, we choose the control at the $n$-th time step  to be an arbitrary $\boldsymbol{\delta}_n^{\Delta t,+}>0$, and for the remaining controls, we choose the elements of a minimizing sequence of controls such that
	\begin{small}
		\begin{align*}
			&   \underset{m \rightarrow \infty}{\lim} \mathbb{E}\left[g^2(\overline{X}_N^{\Delta t})\prod_{k=n+1}^{N-1} L_k^2(\bar{\mathbf{P}}_k,\boldsymbol{\delta}_k^{\Delta t,(m)}(\overline{X}_k^{\Delta t}))\mid \overline{X}_{n+1}^{\Delta t}=\max(0,x +\mathbf{p}^T \boldsymbol{\nu})\right]\\ 
			&= \inf_{\{\boldsymbol{\delta}^{\Delta t}_k\}_{k=n+1,\dots,N-1} \in \mathcal{A}^{N-n-1}} \mathbb{E}\left[g^2(\overline{X}_N^{\Delta t})\prod_{k=n+1}^{N-1} L_k^2(\bar{\mathbf{P}}_k,\boldsymbol{\delta}_k^{\Delta t}(\overline{X}_k^{\Delta t}))\mid \overline{X}_{n+1}^{\Delta t}=\max(0,x +\mathbf{p}^T \boldsymbol{\nu})\right]. 
		\end{align*}
	\end{small}
	
	Therefore, 
	\begin{small}
		\begin{align*}
			&u_{\Delta t}(n,x)\\
			&= \inf_{\{\boldsymbol{\delta}^{\Delta t}_i\}_{i=n,\dots,N-1} \in \mathcal{A}^{N-n}} \Bigg[\exp\left(\left(-2\sum_{j=1}^J a_j(x)+\sum_{j=1}^J\delta_{n,j}^{\Delta t}(x)\right)\Delta t\right)\sum_{\mathbf{p}\in \mathbb{N}^J} \left(\left(\prod_{j=1}^{J} \frac{(\Delta t \cdot \delta_{n,j}^{\Delta t}(x))^{p_j}}{p_j!} (\frac{a_j(x)}{\delta_{n,j}^{\Delta t}(x)})^{2p_j} \right) \right.\\
			& \left. ~~~~~~~~~~~~~~~~~~~~~~~~~~~~~~~~~~~~~~~~~~~~~~ \times  \mathbb{E}\left[g^2(\overline{X}_N^{\Delta t})\prod_{k=n+1}^{N-1} L_k^2(\bar{\mathbf{P}}_k,\boldsymbol{\delta}_k^{\Delta t}(\overline{X}_k^{\Delta t})) \mid \overline{X}_{n+1}^{\Delta t}=\max(0,x +\mathbf{p}^T \boldsymbol{\nu})\right] \right) \Bigg]\\   
			&\leq \exp\left(\left(-2\sum_{j=1}^J a_j(x)+\sum_{j=1}^J\delta_{n,j}^{\Delta t,+}(x)\right)\Delta t\right)\sum_{\mathbf{p}\in \mathbb{N}^J} \left(\left(\prod_{j=1}^{J} \frac{(\Delta t \cdot \delta_{n,j}^{\Delta t,+}(x))^{p_j}}{p_j!} \left(\frac{a_j(x)}{\delta_{n,j}^{\Delta t,+}(x)}\right)^{2p_j} \right) \right.\\
			& \left. ~~~~~~~~~~~~~~~~ \times \inf_{\{\boldsymbol{\delta}^{\Delta t}_k\}_{k=n+1,\dots,N-1} \in \mathcal{A}^{N-n-1}} \mathbb{E}\left[g^2(\overline{X}_N^{\Delta t})\prod_{k=n+1}^{N-1} L_k^2(\bar{\mathbf{P}}_k,\boldsymbol{\delta}_k^{\Delta t}(\overline{X}_k^{\Delta t}))\mid \overline{X}_{n+1}^{\Delta t}=\max(0,x +\mathbf{p}^T \boldsymbol{\nu})\right] \right)
		\end{align*}
		\begin{align*}
			&=\exp\left(\left(-2\sum_{j=1}^J a_j(x)+\sum_{j=1}^J\delta_{n,j}^{\Delta t,+}(x)\right)\Delta t \right)\\
			&~~~~~~~~~~~ \times \sum_{\mathbf{p}\in \mathbb{N}^J}\left(\prod_{j=1}^{J} \frac{(\Delta t \cdot \delta_{n,j}^{\Delta t,+}(x))^{p_j}}{p_j!} \left(\frac{a_j(x)}{\delta_{n,j}^{\Delta t,+}(x)}\right)^{2p_j} \right)\cdot u_{\Delta t}(n+1,\max(0,x+\mathbf{p}^T \boldsymbol{\nu})).
		\end{align*}
	\end{small}
	
	This inequality holds for any arbitrary $\boldsymbol{\delta}_n^{\Delta t,+}>0$, and hence 
	\begin{small}
		\begin{align*}
			u_{\Delta t}(n,x)   \le  \underset{\boldsymbol{\delta}_n^{\Delta t}(x)\in \mathcal{A}_x}{\inf} & \Bigg[\exp\left(\left(-2\sum_{j=1}^J a_j(x)+\sum_{j=1}^J\delta_{n,j}^{\Delta t}(x)\right)\Delta t\right)\\
			&~~~~~~~~~~\sum_{\mathbf{p}\in \mathbb{N}^J}\left(\prod_{j=1}^{J} \frac{(\Delta t \cdot \delta_{n,j}^{\Delta t}(x))^{p_j}}{p_j!} (\frac{a_j(x)}{\delta_{n,j}^{\Delta t}(x)})^{2p_j} \right)\cdot u_{\Delta t}(n+1,\max(0,x+\mathbf{p}^T \boldsymbol{\nu}))\Bigg].
		\end{align*}
	\end{small}
	This completes the proof.
\end{proof}

\section{Proof for Lemma \ref{lem:gradient}}\label{apdx:proofgradient}
The partial derivatives of the second moment $C_{0,\mathbf{x}}\left(\boldsymbol{\delta}^{\Delta t}_n,\dots,\boldsymbol{\delta}^{\Delta t}_{N-1}; \boldsymbol{\beta}\right)$  in \eqref{eq:second_moment} with respect to $\beta_{l}$, $l=1,\dots, (d+1)$ can be expressed as
\begin{small}
	\begin{align}\label{eq:gradientref}
		\frac{\partial}{\partial \beta_l}\mathbb{E}&\left[g^2\left(\overline{\mathbf{X}}_N^{\Delta t,\boldsymbol{\beta}}\right)\prod_{k=0}^{N-1} L_k^2\left(\bar{\mathbf{P}}_k,\hat{\boldsymbol{\delta}}^{\Delta t}(k,\overline{\mathbf{X}}_k^{\Delta t,\boldsymbol{\beta}};\boldsymbol{\beta})\right)\right]\nonumber\\
		&= \frac{\partial}{\partial \beta_l}\mathbb{E}\left[g^2\left(\hat{\mathbf{X}}_N^{\Delta t}\right)\prod_{k=0}^{N-1} L_k\left(\mathbf{P}_k,\hat{\boldsymbol{\delta}}^{\Delta t}(k,\hat{\mathbf{X}}_k^{\Delta t};\boldsymbol{\beta})\right)\right]\nonumber\\
		&\overset{(1)}{=}  \mathbb{E}\left[ \frac{\partial}{\partial \beta_l}\left(g^2\left(\hat{\mathbf{X}}_N^{\Delta t}\right)\prod_{k=0}^{N-1} L_k\left(\mathbf{P}_k,\hat{\boldsymbol{\delta}}^{\Delta t}(k,\hat{\mathbf{X}}_k^{\Delta t};\boldsymbol{\beta})\right)\right)\right]\nonumber\\
		&\overset{(2)}{=} \mathbb{E}\left[ g^2\left(\hat{\mathbf{X}}_N^{\Delta t}\right)\frac{\partial}{\partial \beta_l}\left(\prod_{k=0}^{N-1} L_k\left(\mathbf{P}_k,\hat{\boldsymbol{\delta}}^{\Delta t}(k,\hat{\mathbf{X}}_k^{\Delta t};\boldsymbol{\beta})\right)\right)\right],
	\end{align}
\end{small}
where the Poisson increments with respect to the TL measure are given in $\mathbf{P}_n$ with $(\mathbf{P}_n)_j:=P_{n,j}=\mathcal{P}_{n,j}\left(a_j(\hat{\mathbf{X}}_k^{\Delta t})\Delta t\right)$ for $j=1,\dots,J$.
In $\overset{(1)}{=}$, we assume that the expected value and the derivative commute (see \cite{l1995note} Assumption A1(1) for sufficient conditions). In $\overset{(2)}{=}$, we consider that $g^2\left(\hat{\mathbf{X}}_N^{\Delta t}\right)$ is based on the original TL measure and hence is not dependent on $\beta_l$.  

In \eqref{eq:gradientref}, the term 
$$
\frac{\partial}{\partial \beta_l}\left(\prod_{k=0}^{N-1} L_k\left(\mathbf{P}_k,\hat{\boldsymbol{\delta}}^{\Delta t}(k,\hat{\mathbf{X}}_k^{\Delta t};\boldsymbol{\beta})\right)\right)
$$
is only deterministically dependent on $\boldsymbol{\beta}$, since $\hat{\mathbf{X}}_k^{\Delta t}$ is independent of $\boldsymbol{\beta}$, and $P_{k,j}\sim Poi(a_j(\hat{\mathbf{X}}_k^{\Delta t})\Delta t)$. Thus, the derivative can be computed in a closed form using the identity
\begin{small}
	\begin{align}\label{eq:logderivative}
		\frac{\partial}{\partial x} \ln(f(x))= \frac{1}{f(x)}\frac{\partial}{\partial x} f(x) \iff \frac{\partial}{\partial x} f(x)= f(x)\frac{\partial}{\partial x} \ln(f(x)).
	\end{align}
\end{small}

We compute the derivative from \eqref{eq:logderivative} using the following steps.
\begin{enumerate}
	\item Apply \eqref{eq:logderivative},
	\begin{small}
		\begin{align*}
			\frac{\partial}{\partial \beta_l}&\left(\prod_{k=0}^{N-1} L_k\left(\mathbf{P}_k,\hat{\boldsymbol{\delta}}^{\Delta t}(k,\hat{\mathbf{X}}_k^{\Delta t};\boldsymbol{\beta})\right)\right)\\ 
			&=\left(\prod_{k=0}^{N-1} L_k\left(\mathbf{P}_k,\hat{\boldsymbol{\delta}}^{\Delta t}(k,\hat{\mathbf{X}}_k^{\Delta t};\boldsymbol{\beta})\right) \right)\frac{\partial}{\partial \beta_l} \ln\left( \prod_{k=0}^{N-1} L_k\left(\mathbf{P}_k,\hat{\boldsymbol{\delta}}^{\Delta t}(k,\hat{\mathbf{X}}_k^{\Delta t};\boldsymbol{\beta})\right)\right)\\
			&=\left(\prod_{k=0}^{N-1} L_k\left(\mathbf{P}_k,\hat{\boldsymbol{\delta}}^{\Delta t}(k,\hat{\mathbf{X}}_k^{\Delta t};\boldsymbol{\beta})\right) \right)\sum_{k=0}^{N-1} \frac{\partial}{\partial \beta_l} \ln\left( L_k\left(\mathbf{P}_k,\hat{\boldsymbol{\delta}}^{\Delta t}(k,\hat{\mathbf{X}}_k^{\Delta t};\boldsymbol{\beta})\right)\right) . \\
		\end{align*}
	\end{small}
	\item The remaining derivative can be derived by chain rule, 
	\begin{small}
		\begin{align*}
			\frac{\partial}{\partial \beta_l} \ln\left( L_k\left(\mathbf{P}_k,\hat{\boldsymbol{\delta}}^{\Delta t}(k,\hat{\mathbf{X}}_k^{\Delta t};\boldsymbol{\beta})\right)\right)=\frac{1}{L_k\left(\mathbf{P}_k,\hat{\boldsymbol{\delta}}^{\Delta t}(k,\hat{\mathbf{X}}_k^{\Delta t};\boldsymbol{\beta})\right)} \frac{\partial}{\partial \beta_l}  L_k\left(\mathbf{P}_k,\hat{\boldsymbol{\delta}}^{\Delta t}(k,\hat{\mathbf{X}}_k^{\Delta t};\boldsymbol{\beta})\right).
		\end{align*}
	\end{small}
	\item Apply a second chain rule, 
	\begin{align}\label{eq:deriv}
		\frac{\partial}{\partial \beta_l}  L_k\left(\mathbf{P}_k,\hat{\boldsymbol{\delta}}^{\Delta t}(k,\hat{\mathbf{X}}_k^{\Delta t};\boldsymbol{\beta})\right)=\frac{\partial}{\partial \beta_l} \hat{\boldsymbol{\delta}}^{\Delta t}(k,\hat{\mathbf{X}}_k^{\Delta t};\boldsymbol{\beta}) \cdot \nabla_{\boldsymbol{\delta}} L_k(\mathbf{P}_k,\boldsymbol{\delta}).
	\end{align}
	\item In \eqref{eq:deriv}, we have from \eqref{eq:stepwiselh},
	\begin{align}\label{eq:Lk}
		L_k\left(\mathbf{P}_k,\boldsymbol{\delta}\right)=\exp\left(-\left(\sum_{j=1}^J a_j(\hat{\mathbf{X}}_{k}^{\Delta t})- \delta_j\right)\Delta t\right)  \cdot \prod_{j=1}^J\left(\frac{a_j(\hat{\mathbf{X}}_{k}^{\Delta t})}{\delta_j}\right)^{P_{k,j}},
	\end{align}
	hence
	\begin{small}
		\begin{align*}
			\frac{\partial}{\partial \delta_i} L_k\left(\mathbf{P}_k,\boldsymbol{\delta}\right)&=\Delta t \exp\left(-\left(\sum_{j=1}^J a_j(\hat{\mathbf{X}}_{k}^{\Delta t})- \delta_j\right)\Delta t\right) \cdot \prod_{j=1}^J\left(\frac{a_j(\hat{\mathbf{X}}_{k}^{\Delta t})}{\delta_j}\right)^{P_{k,j}}\\
			&+\exp\left(-\left(\sum_{j=1}^J a_j(\hat{\mathbf{X}}_{k}^{\Delta t})- \delta_j\right)\Delta t\right)  \cdot(-P_{k,i}) \frac{a_i^{P_{k,i}}}{\delta_i^{P_{k,i}+1}}\prod_{j=1,j\neq i}^J\left(\frac{a_j(\hat{\mathbf{X}}_{k}^{\Delta t})}{\delta_j}\right)^{P_{k,j}}\\
			&=\exp\left(-\left(\sum_{j=1}^J a_j(\hat{\mathbf{X}}_{k}^{\Delta t})- \delta_j\right)\Delta t\right)  \prod_{j=1}^J\left(\frac{a_j(\hat{\mathbf{X}}_{k}^{\Delta t})}{\delta_j}\right)^{P_{k,j}}\cdot \left(\Delta t - \frac{P_{k,i}}{\delta_i} \right)\\
			&=L_k(\mathbf{P}_k,\boldsymbol{\delta})\cdot \left(\Delta t - \frac{P_{k,i}}{\delta_i} \right).
		\end{align*}
	\end{small}
	\item In \eqref{eq:deriv}, from \eqref{eq:deltafromu}, we obtain
	\begin{small}
		\begin{align}\label{eq:gradientd}
			\frac{\partial}{\partial \beta_l} \hat{\delta}_j^{\Delta t}(k,\mathbf{x};\boldsymbol{\beta})&=a_j(\mathbf{x})\frac{1}{2}\sqrt{\frac{\hat{u}(\frac{(k+1)\Delta t}{T},\mathbf{x};\boldsymbol{\beta})}{\hat{u}(\frac{(k+1)\Delta t}{T},\max(\mathbf{x}+\nu_j,0)}} \nonumber\\
			&\cdot\left(\frac{\frac{\partial }{\partial \beta_l}\hat{u}(\frac{(k+1)\Delta t}{T},\max(\mathbf{x}+\nu_j,0);\boldsymbol{\beta})}{\hat{u}(\frac{(k+1)\Delta t}{T},\mathbf{x};\boldsymbol{\beta})}\right.\nonumber\\
			&\left.~~~~~~-\frac{\hat{u}(\frac{(k+1)\Delta t}{T},\max(\mathbf{x}+\nu_j,0);\boldsymbol{\beta})\frac{\partial }{\partial \beta_l}\hat{u}(\frac{(k+1)\Delta t}{T},\mathbf{x};\boldsymbol{\beta})}{\hat{u}(\frac{(k+1)\Delta t}{T},\mathbf{x};\boldsymbol{\beta})^2}\right)\nonumber\\
			&=\frac{a_j(\mathbf{x})^2}{2\hat{\delta}_j^{\Delta t}(k,\mathbf{x};\boldsymbol{\beta})}\cdot\left(\frac{\frac{\partial }{\partial \beta_l}\hat{u}(\frac{(k+1)\Delta t}{T},\max(\mathbf{x}+\nu_j,0);\boldsymbol{\beta})}{\hat{u}(\frac{(k+1)\Delta t}{T},\mathbf{x};\boldsymbol{\beta})}\right.\nonumber\\
			&\left.~~~~~~-\frac{\hat{u}(\frac{(k+1)\Delta t}{T},\max(\mathbf{x}+\nu_j,0);\boldsymbol{\beta})\frac{\partial }{\partial \beta_l}\hat{u}(\frac{(k+1)\Delta t}{T},\mathbf{x};\boldsymbol{\beta})}{\hat{u}(\frac{(k+1)\Delta t}{T},\mathbf{x};\boldsymbol{\beta})^2}\right)
		\end{align}
	\end{small}
	where $\frac{\partial }{\partial \beta_l}\hat{u}_{\Delta t}(t,\mathbf{x};\boldsymbol{\beta})$ depends on the chosen ansatz. 
\end{enumerate}

Combining the previous steps, the gradient can be expressed as
\begin{small}
	\begin{align}\label{eq:fullgradient}
		&\frac{\partial}{\partial \beta_l}\left(\prod_{k=0}^{N-1} L_k\left(\mathbf{P}_k,\hat{\boldsymbol{\delta}}^{\Delta t}(k,\hat{\mathbf{X}}_k^{\Delta t};\boldsymbol{\beta})\right)\right)\nonumber\\
		&=\underbrace{\left(\prod_{k=0}^{N-1} L_k\left(\mathbf{P}_k,\hat{\boldsymbol{\delta}}^{\Delta t}(k,\hat{\mathbf{X}}_k^{\Delta t};\boldsymbol{\beta})\right)\right)}_{:=L(\hat{\mathbf{X}}^{\Delta t};\boldsymbol{\beta})}\underbrace{\left(\sum_{k=1}^{N-1}\sum_{j=1}^J \left(\Delta t - \frac{{P}_{k,j}}{\hat{\delta}_j^{\Delta t}(k,\hat{\mathbf{X}}_k^{\Delta t};\boldsymbol{\beta})}\right)\cdot \frac{\partial}{\partial \beta_l} \hat{\delta}_j^{\Delta t}(k,\hat{\mathbf{X}}_k^{\Delta t};\boldsymbol{\beta})\right)}_{:=S(\hat{\mathbf{X}}^{\Delta t};\boldsymbol{\beta})},
	\end{align}
\end{small}
where the gradient of $ \hat{\delta}_j^{\Delta t}$ is dependent on the ansatz used and given by \eqref{eq:gradientd}.

Since the MC estimator \eqref{eq:gradientref} may have a large variance, we again apply IS,
\begin{small}
	\begin{align}\label{eq:ISgradient}
		\frac{\partial}{\partial \beta_l}\mathbb{E}&\left[g^2\left(\overline{\mathbf{X}}_N^{\Delta t,\boldsymbol{\beta}}\right)\prod_{k=0}^{N-1} L_k^2\left(\bar{\mathbf{P}}_k,\hat{\boldsymbol{\delta}}^{\Delta t}(k,\overline{\mathbf{X}}_k^{\Delta t,\boldsymbol{\beta}};\boldsymbol{\beta})\right)\right]\nonumber\\
		\overset{\eqref{eq:gradientref}}{=}
		&\mathbb{E}\left[ g^2\left(\hat{\mathbf{X}}_N^{\Delta t}\right)\frac{\partial}{\partial \beta_l}\left(\prod_{k=0}^{N-1} L_k\left(\mathbf{P}_k,\hat{\boldsymbol{\delta}}^{\Delta t}(k,\hat{\mathbf{X}}_k^{\Delta t};\boldsymbol{\beta})\right)\right)\right]\nonumber\\
		\overset{IS}{=}&\mathbb{E}\left[ \left(\prod_{k=0}^{N-1} L_k\left(\bar{\mathbf{P}}_k,\hat{\boldsymbol{\delta}}^{\Delta t}(k,\overline{\mathbf{X}}_k^{\Delta t,\boldsymbol{\beta}};\boldsymbol{\beta})\right)\right) g^2\left(\overline{\mathbf{X}}_N^{\Delta t,\boldsymbol{\beta}}\right)\frac{\partial}{\partial \beta_l}\left(\prod_{k=0}^{N-1} L_k\left(\bar{\mathbf{P}}_k,\hat{\boldsymbol{\delta}}^{\Delta t}(k,\overline{\mathbf{X}}_k^{\Delta t,\boldsymbol{\beta}};\boldsymbol{\beta})\right)\right)\right]\nonumber\\
		\overset{\eqref{eq:fullgradient}}{=}&\mathbb{E}\left[L(\overline{\mathbf{X}}^{\Delta t,\boldsymbol{\beta}};\boldsymbol{\beta})^2\cdot g^2\left(\overline{\mathbf{X}}_N^{\Delta t,\boldsymbol{\beta}}\right) \cdot S(\overline{\mathbf{X}}^{\Delta t,\boldsymbol{\beta}};\boldsymbol{\beta})\right].
	\end{align}
\end{small}

\end{document}